\documentstyle[12pt,righttag]{amsart}
\newtheorem{theorem}{Theorem\,}[section]
\newtheorem{proposition}{Proposition\,}[section]
\newtheorem{lemma}{Lemma\,}[section]

\newtheorem{remark}{Remark\,}[section]

\begin{document}
\vskip 6cm
\vskip 6cm
\begin{center}
\begin{Large}
Commutation relations of vertex operators related  \\
with the spin representation of $U_q(D_n^{(1)})$
\end{Large}
\end{center}
\vskip 1cm
\begin{center}
Yoshiyuki Koga\footnote{This work was partially supported 
by Research Fellowships of the Japan Society for the Promotion 
of Science for Young Scientists 
and Grant-in-Aid for Scientific Research, the Ministry of Education,
Science and Culture.} 
\end{center}
\vskip 0.5cm
\begin{center}
Department of Mathematics, Graduate School of Science,\\ Osaka University,
Toyonaka, Osaka 560, Japan
\end{center}
\vskip 1.5cm
\begin{abstract}
We calculate commutation relations of vertex operators for the spin representation
of $U_q(D_n^{(1)})$ by using recursive formulae of R-matrices.
In quantum symmetry approach \cite{DFJMN}, 
we obtain the energy and momentum spectrum of 
he quantum spin chain model related with the spin representation 
from these commutation relations.
\end{abstract}

\section{introduction}

A new scheme to study the XXZ spin chain in massive regime
is introduced by Davies et al in \cite{DFJMN}. 
In this scheme, which is called ``quantum symmetry approach", 
the space of states is identified with some $U_q(A_1^{(1)})$-module.
Under this identification, the XXZ Hamiltonian can be constructed on this module
as an operator which commutes with the action of $U_q(A_1^{(1)})$.
Furthermore, the row transfer matrix and creation (annihilation) operators 
can be constructed by using vertex operators.

In a similar way to \cite{DFJMN}, the higher spin chain \cite{IIJMNT} and the models
related to the vector representation of $U_q(A_n^{(1)})$ \cite{DatO}, $U_q(B_n^{(1)})$,
 $U_q(D_n^{(1)})$ \cite{DaiO} are studied. In this paper, we consider the model
related with the {\rm spin} representation of $U_q(D_n^{(1)})$.

In physics the model that we consider here is explained as follows.
It is a one dimensional quantum spin chain model constructed from the spin representation
of $U_q(D_n^{(1)})$. The Hamiltonian acts on the space of the infinite tensor product
\[
\cdots \otimes V^{(n)}_{k+1} \otimes V^{(n)}_k \otimes V^{(n)}_{k-1} \otimes \cdots,
\quad (k \in \Bbb Z),
\]
where $V_{k}^{(n)}$ denotes a copy of spin representation $V^{(+)}$(see Section 2.4). 
The explicit form of the Hamiltonian is given by
\[
{\cal H} = \sum_{k \in Z} h_{k+1,k}, \quad 
h_{k+1,k} = \cdots \otimes \operatorname{id}_{V^{(n)}_{k+2}} \otimes h \otimes
\operatorname{id}_{V^{(n)}_{k-1}} \otimes \cdots,
\]
\[
h = -(q-q^{-1})z \frac{d}{dz} PR(z) |_{z=1},
\]
where $h_{k+1, k}$ acts non-trivially only on the $(k+1)$-th and the $k$-th component.
The operator $R(z)$ is the R-matrix
\[
R(z_1/z_2) \; : \; 
V_{z_1}^{(n)} \otimes V_{z_2}^{(n)} \longrightarrow
V_{z_1}^{(n)} \otimes V_{z_2}^{(n)},
\]
and $P$ denotes the transposition i.e. $P(v \otimes w) = w \otimes v$.

In quantum symmetry approach, the problem is formulated as follows.
Consider the space of states as $U_q(D_n^{(1)})$-module
\begin{equation*}
\underset{\lambda, \mu}{
  \bigoplus}  V(\lambda) \otimes V(\mu)^{*a},
\label{identification}
\end{equation*}
where $\lambda$ and $\mu$ are level one dominant integral weights,
$V(\lambda)$ and $V(\lambda)^{*a}$ denote the irreducible highest weight module with 
highest weight $\lambda$ and its antipode dual (cf. \cite{DFJMN}).
We use vertex operators associated to $U_q(D_n^{(1)})$-modules $V^{(k)}_z$
in order to construct the transfer matrix (the Hamiltonian is obtained from 
the transfer matrix $T(z)$ as follows: ${\cal H}=-(q-q^{-1})zd\log T(z)/dz|_{z=1}$) 
and creation (annihilation) operators.
Here $V^{(k)}_z$ is the affinization of $V^{(k)}$
that is the following representation 
of the ``derived subalgebra" $U'_q(D_n^{(1)})$ (see section 2):
\begin{equation}
 \begin{cases}
   V^{(1)} & \text{the vector representation,}\\ 
   V^{(k)}\;(k = 2,\cdots,n-2) & \text{the fusion representation,}\\
   V^{(n-1)},\; V^{(n)} & \text{the spin representation,}
 \end{cases}
\end{equation}

The row transfer matrix and creation (annihilation) operators are constructed
by using vertex operators for $V^{(k)}_z$ and its dual modules $V^{(k)*a^{\pm 1}}_z$.
Vertex operators for
these modules are defined by the following $U_q(D_n^{(1)})$-homomorphisms:
\begin{eqnarray*}
(\text{type I}) && \Phi^{\mu \; V^{(k)}}_{\lambda}(z) :  
V(\lambda) \longrightarrow V(\mu) \otimes V^{(k)}_z, \\
(\text{type II}) && \Phi^{V^{(k)} \; \mu}_{\lambda}(z) :  
V(\lambda) \longrightarrow V^{(k)}_z \otimes  V(\mu).
\end{eqnarray*}
\begin{eqnarray*}
(\text{type I}) && \Phi^{\mu \; V^{(k)*a^{\pm 1}}}_{\lambda}(z) :  
V(\lambda) \longrightarrow V(\mu) \otimes V^{(k)*a^{\pm 1}}_z, \\
(\text{type II}) && \Phi^{V^{(k)*a^{\pm 1}} \; \mu}_{\lambda}(z) :  
V(\lambda) \longrightarrow V^{(k)*a^{\pm 1}}_z \otimes  V(\mu).
\end{eqnarray*}

Commutation relations of vertex operators give us
commutation relations of the transfer matrix and creation (annihilation) operators,
and then the excitation spectra of the Hamiltonian $\cal H$.
In fact, we can show that vertex operators have the following
commutation relations:
\begin{equation}
  \begin{split}
    \Phi^{\nu \; V^{(n)}}_{\mu'}(z_2) \Phi^{V^{(k)} \; \mu'}_{\lambda}(z_1) 
    &= \tau^{(k)}(z_1/z_2) 
    \Phi^{V^{(k)} \; \nu}_{\mu}(z_1) \Phi^{\mu \; V^{(n)}}_{\lambda}(z_2),\\
    \Phi^{\nu \; V^{(n)}}_{\mu'}(z_2) \Phi^{V^{(k)*a^{-1}} \; \mu'}_{\lambda}(z_1) 
    &= \tau^{(k)}(z_1/z_2)^{-1} 
   \Phi^{V^{(k)*a^{-1}} \; \nu}_{\mu}(z_1) \Phi^{\mu \; V^{(n)}}_{\lambda}(z_2).
  \end{split}
 \label{intro}
\end{equation}
\begin{eqnarray*}
\tau^{(k)}(z) = 
 \begin{cases}
   \;{\displaystyle z^{-\frac{k}{2}} \prod_{j=1}^{k}} 
     \dfrac{\Theta_{q^{4n-4}}(-(-q)^{k+n-2j} z)}{
       \Theta_{q^{4n-4}}(-(-q)^{k+n-2j} z^{-1})}
           \quad & (1 \leq k \leq n-2),\\
   \; & \; \\ 
   \;{\displaystyle z^{-\frac{n-2}{4}} \prod_{j=1}^{[ \frac{n-1}{2} ]} }  
     \dfrac{\Theta_{q^{4n-4}}(-(-q)^{4i-1}z)}{\Theta_{q^{4n-4}}(-(-q)^{4i-1}z^{-1})} 
           \quad & (k=n-1),\\ 
   \; & \; \\
   \; {\displaystyle z^{-\frac{n}{4}} \prod_{j=1}^{[ \frac{n}{2} ]} }  
     \dfrac{\Theta_{q^{4n-4}}(-(-q)^{4i-3}z)}{\Theta_{q^{4n-4}}(-(-q)^{4i-3}z^{-1})}
           \quad & (k=n),
 \end{cases}
\end{eqnarray*}
where $\Theta_p(z)$ is the theta function given by
\[
\Theta_p(z) = (z;p)_{\infty}(p z^{-1};p)_{\infty}(p;p)_{\infty},\quad
 (z;p)_{\infty} = \prod_{i=0}^{\infty}(1-p^i z),
\]
and $[ x ]$ denotes the greatest integer not exceeding $x$.
The explicit form of $\tau^{(k)}(z)$ is crucial to calculate the excitation spectra.
From this formula, we can give explicit forms of the energy $\epsilon^{(k)}(\theta)$
and the momentum $p^{(k)}(\theta)$ with the rapidity variable $\theta$,
\[
e^{-i p^{(k)}(\theta)} = \tau^{(k)}(z), \quad
\epsilon^{(k)}(\theta) = -(q-q^{-1})z\frac{d}{dz}\log\tau^{(k)}(z),\quad -z=e^{2\pi i \theta}.
\] 

To prove the commutation relations \eqref{intro}, 
it is enough to show that the vacuum expectation value of the both sides 
in \eqref{intro} are equal (see section 5.1).
Since the vacuum expectation values satisfy the q-KZ equation,
we can obtain their explicit forms
by finding appropriate solutions of the q-KZ equation (cf. \cite{DFJMN}).

Here we state some comments about excitation spectra. 
By taking the scaling limit as in \cite{BVV} and \cite{DaiO},
we have the following relativistic spectrum:
\begin{equation*}
  \begin{matrix}
  P^{(k)}(\theta) 
    = 2 \mu \sin(\dfrac{\pi k}{2n-2}) \operatorname{sh}(v), 
  &E^{(k)}(\theta) 
    = 2 \mu \sin(\dfrac{\pi k}{2n-2}) \operatorname{ch}(v),
    &\quad (1 \leq k \leq n-2),\\
  P^{(k)}(\theta) = \mu \operatorname{sh}(v), 
  &E^{(k)}(\theta) = \mu  \operatorname{ch}(v),
    &\quad (k = n-1,n),
  \end{matrix}
\end{equation*}
where $P^{(k)}(\theta),\; E^{(k)}(\theta)$ and $v$ can be considered as an appropriate
scaled version of $p^{(k)}(\theta),\; \epsilon^{(k)}(\theta)$ and $\theta$ respectively.
These spectra are exactly same as those of the nonlinear sigma model in \cite{ORW}.
Furthermore, this spectrum coincides with the mass spectrum of 
the spin chain constructed from vector representation of $U_q(D_n^{(1)})$ in \cite{DaiO}.

In general, the structure of the space of states 
turns out to be quite different when we change the region of the parameter
$q$ in the Hamiltonian $\cal H$, 
so we have to discuss the region where we can use the identification \eqref{identification}. 
We are not able to determine this region at this point, and we only state a conjecture.
The region where the identification \eqref{identification} is effective is given by
\[ -1 < q < 0. \] 
These kind of conjectures are already given in \cite{DatO}, \cite{DaiO}. 
Let us consider the case of $n = 3$,
then $U_q(D_3^{(1)})$ is isomorphic to $U_q(A_3^{(1)})$ and the spin representation 
of $U_q(D_3^{(1)})$ is the vector representation of $U_q(A_3^{(1)})$. 
The model related with the vector representation of $U_q(A_n^{(1)})$
is studied by Date and Okado in \cite{DatO}, and our conjecture coincides with theirs
when $n = 3$. 
In this $U_q(A_n^{(1)})$-case, the validity of ``quantum symmetry approach" on this region
is supported by
Bethe Ansatz results \cite{BVV}, however in our case, 
as far as the author knows, there is no similar 
Bethe Ansatz result to support the conjecture.

This paper is organized as follows. We prepare the notation and give 
the definition of the quantum affine algebra $U_q(D_n^{(1)})$ in Section 2. 
In the same section, we construct the spin representations of $U_q(D_n^{(1)})$ and its R-matrices.
By using these R-matrices, the fusion representations are constructed in Section 3.
We also give explicit forms of the R-matrices corresponding to the fusion representations there. 
In Section 4, we give vertex operators for the modules constructed in the previous sections 
and get two point functions.  
Commutation relations of vertex operators are obtained in Section 5. In the same section,
we further make a comment on a relationship between excitation spectra and fusion procedures.
In Section 6, we describe the formulation of the model.
Finally, in Appendix, we give a detailed explanation on calculation of two point functions, 
and construct an isomorphism between the fusion representation and its dual. 

\subsection{Acknowledgment}

The author would like to express his gratitude to Prof. M. Okado,
Prof. K. Nagatomo and Prof. E. Date for their valuable advices. 
He also thanks to Prof. T. Nakanishi for his comment about relation
between fusion and excitation spectra. Furthermore he expresses his thanks
to Prof. M. Ikawa for giving him continuous encouragements. 

\section{Spin representations of the quantum group $U_q(D^{(1)}_n)$}

\subsection{Notation}
In this paper, we use the same notation in \cite{K}.
Let $\frak g$ be the affine Lie algebra of type $D^{(1)}_n$
and $\frak g = {\frak n}_- \oplus {\frak h} \oplus {\frak n}_+$ be the
triangular decomposition.
Let $\alpha_i, \; \alpha_i^{\vee} = h_i, (i = 0, 1, \ldots, n)$
be the simple roots and the simple coroots and $\Lambda_i$ be the fundamental weights
i.e. $\langle \Lambda_i, h_j \rangle = \delta_{i,j}$.
We denote the scaling element and the center by $d$ and $c$
then we can choose elements
$h_1,h_2,\ldots,h_n,c$ and $d$ as a basis of the Cartan subalgebra 
${\frak h}$. Let us define special elements 
$\rho, \delta \in {\frak h}^*$ by $\rho = \sum_{i=0}^n \Lambda_i$ and 
$\delta = \alpha_0 + \alpha_1 + 2\alpha_2 + \cdots + 2\alpha_{n-2} + \alpha_{n-1} + \alpha_n$.
We also denote by $h^{\vee}$ the dual Coxeter number, in $D^{(1)}_n$-case,
this is equal to $2n-2$. 
The invariant bilinear form is normalized by $(\alpha_i|\alpha_i) = 2$.
Let $\overset{\circ}{\frak g}$ be the Lie algebra of type $D_n$, underlying the affine Lie algebra $\frak g$.
For $\lambda \in {\frak h}^*$, we denote by $\bar{\lambda}$ the restriction to
the finite part. We can write $\alpha_i$, $\bar{\Lambda}_i$ and $\bar{\rho}$ by using
orthonormal basis $\{\omega_1, \ldots,\omega_n\}$ of $\bar{\frak h}^*$ (cf. \cite{K}) as follows:
\begin{eqnarray}
 \alpha_i &=& 
  \begin{cases}
    \omega_i - \omega_{i+1} & (1 \leq i \leq n-1),\\
    \omega_{n-1} + \omega_n & (i=n),
  \end{cases}\\
 \bar{\Lambda}_i &=&
  \begin{cases}
    \omega_1+\cdots+\omega_i & (1 \leq i \leq n-2),\\
    \frac{1}{2}(\omega_1+\cdots+\omega_{n-1}-\omega_n) & (i = n-1),\\
    \frac{1}{2}(\omega_1+\cdots+\omega_{n-1}+\omega_n) & (i = n),
  \end{cases}\\
 2\bar{\rho} &=& (2n-2)\omega_1+(2n-4)\omega_2+\cdots+2\omega_{n-1}.
 \label{classical weight}
\end{eqnarray} 
Let us denote $\hat{\sigma_1}$, $\hat{\sigma_2}$ and $\hat{\sigma_3}$ 
the following Dynkin diagram automorphism:
\begin{alignat*}{3}
    \hat{\sigma_1} &:\;\text{$i$-vertex} &\;\longmapsto &\;\text{$(1-i)$-vertex}  \quad &(i&=0,1),  \\
             &:\;\text{$i$-vertex} &\;\longmapsto &\;\text{$i$-vertex} \quad &(i&=2,3,\cdots,n), \\
    \hat{\sigma_2} &:\;\text{$i$-vertex} &\;\longmapsto &\;\text{$i$-vertex} \quad &(i&=0,1,\cdots,n-2), \\
             &:\;\text{$i$-vertex} &\;\longmapsto &\;\text{$(2n-i-1)$-vertex} \quad &(i&=n-1,n), \\
    \hat{\sigma_3} &:\;\text{$i$-vertex} &\;\longmapsto &\;\text{$(n-i)$-vertex} \quad &(i&=0,1,\cdots,n),  
\end{alignat*} 
and further we extend the action of these automorphisms to the weight
lattice $P:= \oplus_{i=1}^n {\Bbb Z} \Lambda_i$ by $\hat{\sigma_k}(\Lambda_i) = \Lambda_{\hat{\sigma_k}(i)}$. 
We also use the following notations:
\begin{align*}
 & [n]_q = \frac{q^n-q^{-n}}{q-q^{-1}},\\
 & {\fracwithdelims[][0pt]{i}{j}}_q 
   = \frac{[i]_q [i-1]_q \cdots [i+1-j]_q}{[j]_q [j-1]_q \cdots [1]_q},\\
 & \xi = q^{h^{\vee}}=q^{2n-2},\\
 & p = q^{2(h^{\vee}+1)} = q^{4n-2}.
\end{align*} 

\subsection{Definition of the quantum group $U_q(D^{(1)}_n)$}
The quantum group $U_q(\frak g)$ is the associative algebra over
${\Bbb Q}(q)$ with generators $e_i$, $f_i$, $t_i=q^{h_i}$ $(i=0,1,\cdots,n)$ and $q^d$.
The defining relations are 
\begin{align*}
 & t_i t_j = t_j t_i,\quad t_i q^d = q^d t_i, \\
 & t_i e_j t_i^{-1} = q^{(\alpha_i|\alpha_j)} e_j,\quad
   t_i f_j t_i^{-1} = q^{-(\alpha_i|\alpha_j)} f_j, \\
 & [e_i,f_j] = \delta_{i,j}\dfrac{t_i-t_i^{-1}}{q-q^{-1}},\\
 & e_i e_j - e_j e_i = f_i f_j - f_j f_i = 0  \hspace{22mm} (\text{if}\;(\alpha_i|\alpha_j)=0),\\
 & 
\begin{aligned}
e_i (e_j)^2 - (q+q^{-1}) e_j e_i e_j + (e_j)^2 e_i = 0 \\
f_i (f_j)^2 - (q+q^{-1}) f_j f_i f_j + (f_j)^2 f_i = 0
\end{aligned}
\quad (\text{if}\;(\alpha_i|\alpha_j)=-1).
\end{align*} 
This algebra $U_q(\frak g)$ has the following Hopf algebra structure:
\begin{align*}
 & \bigtriangleup(e_i) = e_i \otimes 1 + t_i \otimes e_i,\quad
   \bigtriangleup(f_i) = f_i \otimes t_i^{-1} + 1 \otimes f_i, \\
 & \bigtriangleup(t_i) = t_i \otimes t_i,\quad
   \bigtriangleup(q^d) = q^d \otimes q^d, \\
 & a(e_i) = -t_i^{-1} e_i,\; a(f_i)=-f_i t_i, \; a(t_i) =t_i^{-1},\;a(q^d)=q^{-d}.
\end{align*}

We denote by $U_q({\frak n}_+)$ (resp. $U_q({\frak n}_-)$) the subalgebras of $U_q(\frak g)$ 
which are generated by $e_i$ $(i = 0,1,\ldots, n)$ (resp. $f_i$ $(i = 0,1,\ldots, n)$)
and also denote by $U_q'(\frak g)$ the subalgebra 
generated by $e_i$, $f_i$, $t_i$ $(i=0,1,\ldots,n)$.
For any $U_q'(\frak g)$-module $(\pi, V)$, 
we can define $U_q(\frak g)$-module structure $(\pi_z, V_z)$ by
\begin{equation}
  \begin{split}
   &V_z = {\Bbb Q}(q)[z,z^{-1}] \otimes V,\\
   &\quad \pi_z(e_i) (z^n \otimes v) = z^{n+\delta_{i,0}} \otimes \pi(e_i) v,\\
   &\quad \pi_z(f_i) (z^n \otimes v) = z^{n-\delta_{i,0}} \otimes \pi(f_i) v,\\
   &\quad \pi_z(t_i) (z^n \otimes v) = z^n \otimes \pi(t_i) v,\\
   &\quad \pi_z(q^d) (z^n \otimes v) = (q z)^n \otimes v,
  \end{split}
  \label{affinization}
\end{equation}
and this $U_q(\frak g)$-module $(\pi_z, V_z)$ is called the affinization of $(\pi, V)$.

The Dynkin diagram automorphisms can be extended to algebra 
automorphisms of $U_q'(\frak g)$ as follows:
\begin{alignat*}{4}
     \hat{\sigma_1}(x_i) &= x_{1-i} & (i&=0,1), & \hat{\sigma_1}(x_i)&=x_i & (i&=2,3,\cdots,n),\\
     \hat{\sigma_2}(x_i) &= x_i & (i&=0,1,\cdots,n-2), \;\;& \hat{\sigma_2}(x_i)&=x_{2n-i-1} & (i&=n-1,n),\\
     \hat{\sigma_3}(x_i) &= x_{n-i} & (i&=0,1,\cdots,n),& & & &
\end{alignat*}
where $x_i$ stands for $e_i$, $f_i$ or $t_i$.

\subsection{R-matrices}
Let $(\pi^{W_1},W_1)$ and $(\pi^{W_2},W_2)$ be finite-dimensional $U_q'(\frak g)$-modules.
An operator $R(z_1/z_2) \in \operatorname{End}(W_1 \otimes W_2)$ is called R-matrix
if it has the following intertwining property:
\begin{equation}
  \begin{split}
     &R^{W_1,W_2}(z_1/z_2)(\pi^{W_1}_{z_1} \otimes \pi^{W_2}_{z_2}) \bigtriangleup(x) \\
     & \qquad = (\pi^{W_1}_{z_1} \otimes \pi^{W_2}_{z_2}) \bigtriangleup'(x) R^{W_1,W_2}(z_1/z_2)
       \quad (x \in U_q(\frak g) ),
   \end{split}
   \label{R-matrix}
\end{equation}
where $\bigtriangleup' = P \circ \bigtriangleup$.

Let $(\pi^{W_3},W_3)$ be the third representation of $U_q'(\frak g)$. The R-matrices satisfy
the Yang-Baxter equation on $W_1 \otimes W_2 \otimes W_3$,
\begin{equation}
  \begin{split}
     &R^{W_1,W_2}(z_1/z_2) R^{W_1,W_3}(z_1/z_3) R^{W_2,W_3}(z_2/z_3) \\
     & \qquad= R^{W_2,W_3}(z_2/z_3) R^{W_1,W_3}(z_1/z_3) R^{W_1,W_2}(z_1/z_2),
  \end{split}
  \label{YB-equation}
\end{equation}
where $R^{W_i,W_j}$ acts non-trivially on the $i$-th and the $j$-th components.

The modified universal R-matrix ${\cal R}'(z)$ is defined by
\begin{equation}
 {\cal R}'(z) 
   = q^{-\sum_{i=1}^n h_i \otimes {\bar{\Lambda}}_i} \sum_{\beta,i} z^{<d,\beta>}
   u_{\beta,i} \otimes u_{-\beta}^i \;(\in U'_q(\frak g) \hat{\otimes} U'_q(\frak g)), 
 \label{modified R-matrix}
\end{equation}
where $\{u_{\beta,i}\}$ (resp. $\{u_{-\beta}^i\}$) are dual bases of weight $\beta$ 
(resp. $-\beta$) component of the subalgebra $U_q({\frak n}_+)$ (resp. $U_q({\frak n}_-)$) and $\hat{\otimes}$ 
means completion in formal topology with respect to the weight decomposition of $U_q({\frak n}_+)$ (cf. \cite{IIJMNT}).
In the same way as \cite{IIJMNT}, we can determine the image of 
the modified universal R matrix under the representation
\[
\pi^{W_1} \otimes \pi^{W_2}\;:\;
U_q'(\frak g) \hat{\otimes} U_q'(\frak g) \longrightarrow
\operatorname{End}(W_1 \otimes W_2).
\]
If our R-matrix on $W_1 \otimes W_2$ is uniquely determined up to multiple scalar factor
then the image of ${\cal R}'(z)$ has the following form: 
\begin{equation}
  (\pi^{W_1} \otimes \pi^{W_2})({\cal R}'(z))
  = {\beta}^{W_1,W_2}(z) R^{W_1,W_2}(z). 
  \label{image of mR-matrix}
\end{equation}
To calculate the scalar factor ${\beta}^{W_1,W_2}(z)$, we have to obtain
the explicit form of a scalar function ${\alpha}^{W_1,W_2}(z)$ appeared 
in the second inversion relation:
\begin{equation}
  {\alpha}^{W_1,W_2}(z) (((R^{W_1,W_2}(z)^{-1})^{t_1})^{-1})^{t_1}
   = (q^{-2\rho}\otimes 1) R^{W_1,W_2}({\xi}^{-2}z) (q^{2\rho}\otimes 1),
  \label{second inversion}
\end{equation}
where the symbol ${t_1}$ means the transpose in the first component.
It is known that the scalar function ${\alpha}^{W_1,W_2}(z)$ for the image of ${\cal R}'(z)$
is equal to one i.e.
\[
( ( ( ({\beta}^{W_1,W_2}(z) R^{W_1,W_2}(z))^{-1} )^{t_1})^{-1})^{t_1}
 = (q^{-2\rho}\otimes 1) {\beta}^{W_1,W_2}({\xi}^{-2}z) R^{W_1,W_2}({\xi}^{-2}z) (q^{2\rho}\otimes 1).
\]
Comparing this equation with \eqref{second inversion}, we have
\[
{\alpha}^{W_1,W_2}(z) {\beta}^{W_1,W_2}(z)^{-1} = {\beta}^{W_1,W_2}({\xi}^{-2}z)^{-1}.
\]
We can determine the explicit form of ${\beta}^{W_1,W_2}(z)$ as a solution of this difference equation.
In the case related with the spin representation, more detailed explanation
on calculation of ${\alpha}^{W_1,W_2}(z)$ and ${\beta}^{W_1,W_2}(z)$ is given in Appendix.

\subsection{Spin representations}
Let $V_{1/2}$ be a 2-dimensional vector space 
spanned by vectors $v_{1/2}$ and $v_{-1/2}$ over the field ${\Bbb Q}(q^{1/2})$.
We define operators $X^+$, $X^-$ and $T$ acting on $V_{1/2}$ by
\[
X^+ v_{\gamma} = v_{\gamma+1}, \; X^- v_{\gamma} = v_{\gamma-1}, \;
Tv_{\gamma} = q^{\gamma}v_{\gamma},
\]
where if ${\gamma} \not= \pm 1/2$ then $v_{\gamma}=0$ i.e. matrices corresponding to $X^+$, $X^-$ and $T$ are
given by
\begin{equation*}
  \begin{pmatrix}
     0&1\\
     0&0
  \end{pmatrix},
  \quad
  \begin{pmatrix}
     0&0\\
     1&0
  \end{pmatrix},
  \quad
  \begin{pmatrix}
     q^{\frac{1}{2}}&0\\
     0&q^{-\frac{1}{2}}
  \end{pmatrix}. 
\end{equation*}
By using these operators,
we can define a representation of $U_q'(\frak g)$ on a vector space 
$V^{(sp)} = (V_{1/2})^{\otimes n}$ as follows \cite{O}: 
\begin{equation}
  \begin{split}
    &\pi^{(sp)}(e_0) 
      = X^- \otimes X^- \otimes 1 \otimes \cdots \otimes 1, \\
    &\pi^{(sp)}(t_0) 
      = T^{-1} \otimes T^{-1} \otimes 1 \otimes \cdots \otimes 1, \\
    &\pi^{(sp)}(e_i) 
      = 1 \otimes \cdots \otimes 1 \otimes \overset{\text{$i$-th}}{X^+} \otimes \overset{\text{$(i+1)$-th}}{X^-} 
          \otimes 1 \otimes \cdots \otimes 1 \;\;(1 \leq i \leq n-1), \\
    &\pi^{(sp)}(t_i) 
      = 1 \otimes \cdots \otimes 1 \otimes \overset{\text{$i$-th}}{T} \otimes \overset{\text{$(i+1)$-th}}{T^{-1}}  
          \otimes 1 \otimes \cdots \otimes 1 \quad (1 \leq i \leq n-1),\\
    &\pi^{(sp)}(e_n) 
      = 1 \otimes \cdot \cdot \cdot \otimes 1 \otimes X^+ \otimes X^+,\\
    &\pi^{(sp)}(t_n) 
      = 1 \otimes \cdot \cdot \cdot \otimes 1 \otimes T \otimes T ,\\
    &\pi^{(sp)}(f_i) = \pi^{(sp)}(e_i)^t. 
  \end{split}
  \label{spin-repre}
\end{equation}
We denote a vector $v_{\gamma_1} \otimes \cdots \otimes v_{\gamma_n} \in V^{(sp)}$ 
by $v_{(\varepsilon_1, ... ,\varepsilon_n)}$ ($\varepsilon_k = \operatorname{sgn}(\gamma_k)$).
The weight of the element $v_{(\varepsilon_1, ... ,\varepsilon_n)}$ is given by 
$\sum_{i=1}^n \varepsilon_i \omega_i/2$, where $\{\omega_k\}$ is the orthogonal basis of $\bar{\frak h}^*$.
Let $V^{(+)}$ and $V^{(-)}$ be the subspaces spanned by
$\{ v_{(\varepsilon_1, ... ,\varepsilon_n)} | \prod_{i=1}^n \varepsilon_i = + \}$ and 
$\{ v_{(\varepsilon_1, ... ,\varepsilon_n)} | \prod_{i=1}^n \varepsilon_i = - \}$,
then we can easily show that they are irreducible submodules. 
These irreducible representations are called spin representations and are denoted by $(\pi^{(\pm)},V^{(\pm)})$.

There exists isomorphisms of 
$U_q(\frak g)$-module (cf. \cite{DaiO})
\begin{equation}
  C^{(sp)}_{\pm} \; : \; V^{(sp)}_{z \xi^{\mp}} \longrightarrow (V^{(sp)}_z)^{*a^{\pm1}}.
  \label{spin-dual}
\end{equation}
We denote the restrictions of $C^{(sp)}_{\pm}$ to the irreducible 
components $V^{(+)}$ (resp. $V^{(-)}$) by 
$C^{(+)}_{\pm}$ (resp. $C^{(-)}_{\pm}$) i.e.
\begin{equation}
 \begin{split}
   &C^{(+)}_{\pm} \; : \; V^{(+)}_{z \xi^{\mp}} 
    \longrightarrow (V^{(\varepsilon)}_z)^{*a^{\pm1}},\\
   &C^{(-)}_{\pm} \; : \; V^{(-)}_{z \xi^{\mp}} 
    \longrightarrow (V^{(-\varepsilon)}_z)^{*a^{\pm1}},
 \end{split}
 \label{dual for spin}
\end{equation}
where $\varepsilon$ is given by
\begin{equation*}
\varepsilon =
 \begin{cases}
   + & (n : \text{even}),\\
   - & (n : \text{odd}).
 \end{cases}
\end{equation*}

We can define the action of the Dynkin diagram automorphisms 
$\hat{\sigma_1}$, $\hat{\sigma_2}$ and $\hat{\sigma_3}$ to the spin representations 
as follows:
\begin{equation}
  \begin{split}
   & \hat{\sigma_1}(v_{(\varepsilon_1,\varepsilon_2,\cdots,\varepsilon_n)})
       = v_{(-\varepsilon_1,\varepsilon_2,\cdots,\varepsilon_n) },\\
   & \hat{\sigma_2}(v_{(\varepsilon_1,\cdots,\varepsilon_{n-1},\varepsilon_n)})
       = v_{(\varepsilon_1,\cdots,\varepsilon_{n-1},-\varepsilon_n)},\\
   & \hat{\sigma_3}(v_{(\varepsilon_1,\varepsilon_2,\cdots,\varepsilon_n)})
       = v_{(-\varepsilon_n,\cdots,-\varepsilon_2,-\varepsilon_1)}.
  \end{split}
\end{equation} 

\subsection{R-matrices related with the spin representations}
In this subsection, we describe recursive forms of the R-matrices related with
the spin representations.   In \cite{O}, the R-matrices related with
the spin representations of $U_q(B_n^{(1)})$ and $U_q(D_n^{(1)})$ are determined,
and the R-matrix for $U_q(B_{n}^{(1)})$ is recursively expressed by using 
the one for $U_q(B_{n-1}^{(1)})$.
We give a similar recursive relation in the case of $U_q(D_n^{(1)})$.

We normalize the R-matrix $R^{V^{(\varepsilon_1)},V^{(\varepsilon_2)}}(z)$ by
\[
R^{V^{(\varepsilon_1)},V^{(\varepsilon_2)}}(z)
v_{(+,+,...,+,\varepsilon_1)} \otimes v_{(+,+,...,+,\varepsilon_2)} =
v_{(+,+,...,+,\varepsilon_1)} \otimes v_{(+,+,...,+,\varepsilon_2)},
\]
and denote this R-matrix by $\bar{R}^{(\varepsilon_1,\varepsilon_2)}_n(z)$
where $n$ stands for the rank of $\bar{\frak g}$.
Expressing the rank of $\frak g$ clearly, we denote the spin representations $V^{(\pm)}$
by $V^{(\pm)}_n$ only in this subsection.
In order to describe recursive formulae of R-matrices, we consider an isomorphism of vector spaces
\begin{equation}
V^{(\varepsilon_1)}_n \otimes V^{(\varepsilon_2)}_n \longrightarrow 
\sum
  \begin{Sb}
     \eta_1,\eta_2=\pm \\
     \eta_1' = \varepsilon_1 \eta_1 \\ 
     \eta_2' = \varepsilon_2 \eta_2
  \end{Sb}
    {\Bbb Q}(q) v_{\eta_1 , \eta_2} \otimes 
    (V^{(\eta_1')}_{n-1} \otimes V^{(\eta_2')}_{n-1}),
  \label{new_expression}
\end{equation}
which maps 
\[
v_{(\nu_1,\nu_2 ... ,\nu_n)} \otimes v_{(\mu_1,\mu_2 ... ,\mu_n)} \; \mapsto \; 
{v_{\nu_1,\mu_1}} \otimes {(v_{(\nu_2 ... ,\nu_n)} \otimes v_{(\mu_2 ... ,\mu_n)})},
\]
where we denote $v_{\nu_1} \otimes v_{\mu_1}$ by $v_{\nu_1,\mu_1}$.
We put
\begin{equation}
a(z) = \frac{q(1-z)}{1-q^2z}, \quad b(z) = \frac{1-q^2}{1-q^2z},\quad
c(z) = \frac{b(z)}{a(z)}.
\label{abc}
\end{equation}
By using this isomorphism, we have the following expressions of R-matrices:

\vspace{5mm}
\noindent for $ n \geq 2 $
\begin{equation}
  \begin{split}
  &\bar{R}^{(\varepsilon,\varepsilon)}_n(z) v_{++} \otimes u_{n-1}
      = v_{++} \otimes \bar{R}^{(\varepsilon,\varepsilon)}_{n-1}(z) u_{n-1},\\
  &\bar{R}^{(\varepsilon,\varepsilon)}_n(z) v_{+-} \otimes u_{n-1}\\
  &\qquad = v_{+-} \otimes a(z)\bar{R}^{(\varepsilon,-\varepsilon)}_{n-1}(q^2z) u_{n-1}
      +v_{-+} \otimes zb(z) \bar{R}^{(-\varepsilon,\varepsilon)}_{n-1}(q^2z) 
     {\sigma}_{n-1} X^{(\varepsilon,-\varepsilon)}_{n-1} u_{n-1}, \\
  &\bar{R}^{(\varepsilon,\varepsilon)}_n(z) v_{-+} \otimes u_{n-1}\\
  &\qquad =v_{+-} \otimes b(z) \bar{R}^{(\varepsilon,-\varepsilon)}_{n-1}(q^2z) 
     {\sigma}_{n-1} X^{(-\varepsilon,\varepsilon)}_{n-1} u_{n-1} 
     +v_{-+} \otimes a(z)\bar{R}^{(-\varepsilon,\varepsilon)}_{n-1}(q^2z) u_{n-1},\\
  &\bar{R}^{(\varepsilon,\varepsilon)}_n(z) v_{--} \otimes u_{n-1}
    = v_{--} \otimes \bar{R}^{(-\varepsilon,-\varepsilon)}_{n-1}(z) u_{n-1},
  \end{split}
  \label{recursive form ++}
\end{equation}
\begin{equation}
  \begin{split}
   &\hspace{-10mm}\bar{R}^{(\varepsilon,-\varepsilon)}_n(z) v_{++} \otimes u_{n-1}
      = v_{++} \otimes \bar{R}^{(\varepsilon,-\varepsilon)}_{n-1}(z) u_{n-1},\\
   &\hspace{-10mm}\bar{R}^{(\varepsilon,-\varepsilon)}_n(z) v_{+-} \otimes u_{n-1}\\
   &\hspace{-10mm}\qquad = v_{+-} \otimes \bar{R}^{(\varepsilon,\varepsilon)}_{n-1}(q^2z) u_{n-1}
      +v_{-+} \otimes zc(z) \bar{R}^{(-\varepsilon,-\varepsilon)}_{n-1}(q^2z) 
     {\sigma}_{n-1} X^{(\varepsilon,\varepsilon)}_{n-1} u_{n-1},\\
   &\hspace{-10mm}\bar{R}^{(\varepsilon,-\varepsilon)}_n(z) v_{-+} \otimes u_{n-1}\\
   &\hspace{-10mm}\qquad =v_{+-} \otimes c(z) \bar{R}^{(\varepsilon,\varepsilon)}_{n-1}(q^2z) 
     {\sigma}_{n-1} X^{(-\varepsilon,-\varepsilon)}_{n-1} u_{n-1} 
     +v_{-+} \otimes \bar{R}^{(-\varepsilon,-\varepsilon)}_{n-1}(q^2z) u_{n-1},\\
   &\hspace{-10mm}\bar{R}^{(\varepsilon,-\varepsilon)}_n(z) v_{--} \otimes u_{n-1}
     = v_{--} \otimes \bar{R}^{(-\varepsilon,\varepsilon)}_{n-1}(z) u_{n-1},
  \end{split}
  \label{recursive form +-}
\end{equation}

\noindent for $ n = 1 $
\begin{eqnarray*}
    \bar{R}^{(\varepsilon,\varepsilon)}_1(z) = 1,  \quad
    \bar{R}^{(\varepsilon,-\varepsilon)}_1(z) =1,
\end{eqnarray*}
where an involution ${\sigma}_{n}$ is defined by
\[
{\sigma}_{n}( v_{(\varepsilon_1,\cdots,\varepsilon_n)} \otimes v_{(\eta_1,\cdots,\eta_n)} ) 
= v_{(\varepsilon_1,\cdots,\varepsilon_{n-1},-\varepsilon_n)} 
\otimes v_{(\eta_1,\cdots,\eta_{n-1},-\eta_n)},
\]
and $X^{(\varepsilon_1,\varepsilon_2)}_n$ is a linear operator on
 $V^{(\varepsilon_1)}_n \otimes V^{(\varepsilon_2)}_n$ that is defined inductively 
as follows:

\vspace{5mm}
\noindent for $ n \geq 2 $
\begin{equation}
  \begin{split}
  & X^{(\varepsilon_1,\varepsilon_2)}_n v_{++} \otimes u_{n-1}
    = v_{++} \otimes X^{(\varepsilon_1,\varepsilon_2)}_{n-1} u_{n-1} \\
  & X^{(\varepsilon_1,\varepsilon_2)}_n v_{+-} \otimes u_{n-1}
    = -q^{-1} v_{+-} \otimes X^{(\varepsilon_1,-\varepsilon_2)}_{n-1} u_{n-1}
      + v_{-+} \otimes {\sigma}_{n-1} u_{n-1},\\
  & X^{(\varepsilon_1,\varepsilon_2)}_n v_{-+} \otimes u_{n-1}
    = v_{+-} \otimes {\sigma}_{n-1} u_{n-1}
      - q v_{-+} \otimes X^{(-\varepsilon_1,\varepsilon_2)}_{n-1} u_{n-1}, \\
  & X^{(\varepsilon_1,\varepsilon_2)}_n v_{--} \otimes u_{n-1}
    = v_{--} \otimes X^{(-\varepsilon_1,-\varepsilon_2)}_{n-1} u_{n-1},
  \end{split}
  \label{def of X}
\end{equation}

\noindent for $ n = 1 $
\begin{eqnarray*}
  X^{(\varepsilon,\varepsilon)}_1 = 0,& & 
  X^{(\varepsilon,-\varepsilon)}_1 =1.
\end{eqnarray*}

We denote by ${\alpha}^{(\varepsilon_1,\varepsilon_2)}_n(z)$ 
the scalar function $\alpha^{V^{(\varepsilon_1)},V^{(\varepsilon_2)}}(z)$
in \eqref{second inversion}.
Then we can find
\begin{equation}
{\alpha}^{(\varepsilon_1,\varepsilon_2)}_n(z) = 
 \begin{cases}
 {\displaystyle \prod_{i=1}^{[ \frac{n}{2} ]} }
 \dfrac{(1-q^{-4i+4}z)(1-q^{-4n+4i}z)}{(1-q^{-4i+2}z)(1-q^{-4n+4i+2}z)}
 \quad & (\varepsilon_1 = \varepsilon_2),\\
 \; & \; \\
 {\displaystyle \prod_{i=1}^{[ \frac{n-1}{2} ]} }
 \dfrac{(1-q^{-4i+2}z)(1-q^{-4n+4i+2}z)}{(1-q^{-4i}z)(1-q^{-4n+4i+4}z)}
 \quad & (\varepsilon_1 \neq \varepsilon_2).
 \end{cases}
 \label{second_inversion_spin}
\end{equation}
This fact is proved in Appendix A.

We also denote by ${\beta}^{(\varepsilon_1,\varepsilon_2)}_n(z)$ 
the scalar factor ${\beta}^{V^{(\varepsilon_1)},V^{(\varepsilon_2)}}(z)$ in \eqref{image of mR-matrix}.
By using the explicit form of the function ${\alpha}^{(\varepsilon_1,\varepsilon_2)}_n(z)$, 
we can show (see Appendix B)
\begin{equation}
  {\beta}^{(\varepsilon_1,\varepsilon_2)}_n(z) =
   \begin{cases}
   {\displaystyle q^{-\frac{n}{4}}\prod_{i=1}^{[ \frac{n}{2} ]} }
      \dfrac{(q^{4n-4i-2}z;{\xi}^2)_{\infty}(q^{4i-2}z;{\xi}^2)_{\infty}}{
      (q^{4n-4i}z;{\xi}^2)_{\infty}(q^{4i-4}z;{\xi}^2)_{\infty}}
      \quad & (\varepsilon_1 = \varepsilon_2),\\
    \; & \; \\
   {\displaystyle  q^{-\frac{n}{4}+\frac{1}{2}}\prod_{i=1}^{[ \frac{n-1}{2} ]} }
      \dfrac{(q^{4n-4i-4}z;{\xi}^2)_{\infty}(q^{4i}z;{\xi}^2)_{\infty}}{
      (q^{4n-4i-2}z;{\xi}^2)_{\infty}(q^{4i-2}z;{\xi}^2)_{\infty}}
      \quad & (\varepsilon_1 \neq \varepsilon_2).
   \end{cases}
   \label{image_MR_spin}
\end{equation}

\section{Fusion construction}

We construct the fusion representations and R-matrices related with them.

\subsection{Construction of $V^{(k)}$}
Since $V^{(+)}$ and $V^{(-)}$ are the highest weight
modules corresponding to the highest weights $\bar{\Lambda}_n$ and $\bar{\Lambda}_{n-1}$ as 
$U_q(\overset{\circ}{\frak g})$-modules, we denote $V^{(+)}$ and $V^{(-)}$ by $V^{(n)}$ and $V^{(n-1)}$. 
We also denote $\bar{R}^{(+,+)}_n(z)$ and $\bar{R}^{(+,-)}_n(z)$
by $\bar{R}^{(n,n)}(z)$ and $\bar{R}^{(n,n-1)}(z)$ respectively.

The R-matrix $\bar{R}^{(n,n)}(z)$ (resp. $\bar{R}^{(n,n-1)}(z)$) have a pole of order one 
at  $z= q^{-2n+2k +2}$ $(k \equiv n \pmod{2})$ (resp. $z = q^{-2n+2k +2}$ 
$(k \not\equiv n \pmod{2})$\;) for all $k=1,2,\ldots,n-2$. We eliminate these poles
by multiplying an appropriate scalar factor to the R-matrix as follows.
\begin{equation*}
  \begin{split}
  &R^{(n,n)}(z) 
    = \prod_{i=1}^{[ \frac{n}{2} ]}(1-q^{4i-2}z)
      \bar{R}^{(n,n)}(z),\\
  &R^{(n,n-1)}(z) 
    = \prod_{i=1}^{[ \frac{n-1}{2} ]}(1-q^{4i}z)
      \bar{R}^{(n,n-1)}(z).
  \end{split}
\end{equation*}
By using these R-matrices,
let us define operators $T^{(k)} \; ( 1\leq k \leq n-2)$ by 
\begin{eqnarray*}
  T^{(k)} = 
    \begin{cases}
    \;\; R^{(n,n)}(q^{-2n+2k+2}) \quad \in \operatorname{End}( V^{(n)} \otimes V^{(n)} )
                                     \quad & ( \; k \equiv n \; \pmod{2}), \\
    \;\; R^{(n,n-1)}(q^{-2n+2k+2}) \in \operatorname{End}( V^{(n)} \otimes V^{(n-1)} )
                                     \quad & ( \; k \not\equiv n \; \pmod{2}).
    \end{cases} 
\end{eqnarray*} 
The following function $\varphi^{(n)}:{\Bbb N}\rightarrow\{n-1,n\}$ is useful in this section:
\begin{equation}
  \varphi^{(n)}(i) = 
  \begin{cases}
    n & (i:\text{even}), \\
    n-1 &   (i:\text{odd}).
  \end{cases}
  \label{numerical function}
\end{equation}
We introduce a vector space
\[
V^{(k)} = (V^{(n)} \otimes V^{(n')})/\ker T^{(k)},
\]
where $1\leq k \leq n-2$ and $n' = \varphi^{(n)}(n-k)$. 
We define an action of $U'_q(\frak g)$
\[
\tilde{\pi}_z^{(k)}\; : \; U'_q(\frak g) 
   \longrightarrow \operatorname{End}( V^{(n)} \otimes V^{(n')} ),
\]
by
\[
\tilde{\pi}_z^{(k)}(x) 
= (\pi_{(-q)^{-n+k+1}z}^{(n)} \otimes \pi_{(-q)^{n-k-1}z}^{(n')})
      \circ \bigtriangleup(x),
\]
then we have the following proposition.

\begin{proposition}
For all $x \in U'_q(\frak g)$,
\begin{equation}
  \tilde{\pi}_z^{(k)}(x) \ker T^{(k)} \subset \ker T^{(k)}.
  \label{action of fusion repre}
\end{equation}
Then $\tilde{\pi}_z^{(k)}$ induces a representation of $U'_q(\frak g)$ on $V^{(k)}$.
\end{proposition}
\noindent
Proof. \ Let $v$ be an arbitrary element in $\ker T^{(k)}$.
By means of the intertwiner property of the R-matrices in \eqref{R-matrix}, we have
\begin{eqnarray*}
T^{(k)} \tilde{\pi}_z^{(k)}(x) v
&=& R^{(n,n')}(q^{-2n+2k+2})
    (\pi_{(-q)^{-n+k+1}z}^{(n)} \otimes \pi_{(-q)^{n-k-1}z}^{(n')})
      \bigtriangleup(x) v \\
&=& (\pi_{(-q)^{-n+k+1}z}^{(n)} \otimes \pi_{(-q)^{n-k-1}z}^{(n')})
      \bigtriangleup'(x)
    R^{(n,n')}(q^{-2n+2k+2}) v \\
&=&0,
\end{eqnarray*}
for all $x \in U'_q(\frak g)$.
Then $\tilde{\pi}_z^{(k)}(x) v \in \ker T^{(k)}$.
\hfill Q.E.D
\vspace{5mm}

We define an action of $q^d$ on $V^{(k)}$ which is 
a representation of $U'_q(\frak g)$ induced by the above $\tilde{\pi}^{(k)}_z$
as similar to \eqref{affinization} and denote it by $(\pi^{(k)}_z,V^{(k)}_z)$.
The following isomorphisms of $U_q(\frak g)$-module are known in \cite{DaiO}:
\begin{equation}
  C^{(k)}_{\pm} \;:\; V^{(k)}_{z\xi^{\mp}} \longrightarrow (V^{(k)}_z)^{*a^{\pm 1}}
    \quad (k = 1,2,\cdots,n-2).
  \label{dual for fusion}
\end{equation}
Explicit forms of $C^{(k)}_{\pm}$ are given in Appendix.

\subsection{R-matrices related with $V^{(k)}$}
For $k = 1,2,\ldots,n-2$, we explicitly construct R-matrices on 
$V^{(k)}\otimes V^{(n)}$, $V^{(k)}\otimes V^{(n-1)}$, $V^{(n)} \otimes V^{(k)}$
and $V^{(n-1)} \otimes V^{(k)}$ (see \eqref{R-matrix}).
Let $m$ be $n$ or $n-1$. We define operators 
\begin{equation*}
  \begin{split}
     &R^{(m,k)}(z) \in
      \operatorname{End}(V^{(m)} \otimes V^{(n)} \otimes V^{(n')}), \\
     &R^{(k,m)}(z) \in
     \operatorname{End}(V^{(n)} \otimes V^{(n')} \otimes V^{(m)}), 
  \end{split}
\end{equation*}
by
\begin{equation*}
  \begin{split}
    &R^{(m,k)}(z) =
    {\bar{R}}^{(m,n')}((-q)^{-n+k+1}z)_{13}
    {\bar{R}}^{(m,n)}((-q)^{n-k-1}z)_{12}, \\
    &R^{(k,m)}(z) =
    {\bar{R}}^{(n,m)}((-q)^{-n+k+1}z)_{13}
    {\bar{R}}^{(n',m)}((-q)^{n-k-1}z)_{23}, \\
  \end{split}
\end{equation*}
where  $n' = \varphi^{(n)}(n-k)$ 
and the subscripts of these R-matrices indicate 
the components that each operator acts on, that
is, $R_{ij}$ non-trivially acts on the $i$-th and the $j$-th components. 

Here we will show that 
\[
R^{(m,k)}(z) ( V^{(m)} \otimes \ker T^{(k)}) \subset V^{(m)} \otimes \ker T^{(k)},
\]
and so $R^{(m,k)}(z)$ (resp. $R^{(k,m)}(z)$) define
operators on $V^{(m)} \otimes V^{(k)}$ (resp. $V^{(k)} \otimes V^{(m)}$).
In fact, for arbitrary $v \in V^{(m)} \otimes \ker T^{(k)}$,
by using the Yang-Baxter equation \eqref{YB-equation}, we have
\begin{eqnarray*}
&&(\operatorname{id} \otimes T^{(k)}) R^{(m,k)}(z) v\\
&&\quad =R^{(n,n')}(q^{-2n+2k+2})_{23} 
   {\bar{R}}^{(m,n')}((-q)^{-n+k+1}z)_{13}
   {\bar{R}}^{(m,n)}((-q)^{n-k-1}z)_{12} v\\
&&\quad ={\bar{R}}^{(m,n)}((-q)^{n-k-1}z)_{12}
   {\bar{R}}^{(m,n')}((-q)^{-n+k+1}z)_{13}
   R^{(n,n')}(q^{-2n+2k+2})_{23} v\\
&&\quad =0.
\end{eqnarray*}
Hence  $(\operatorname{id} \otimes T^{(k)}) R^{(m,k)}(z) v \in V^{(m)} \otimes \ker T^{(k)}$.
We can prove 
\[
R^{(k,m)}(z) ( \ker T^{(k)} \otimes V^{(m)}) \subset \ker T^{(k)} \otimes V^{(m)}
\]
similarly. Then we have
\begin{proposition}
$R^{(m,k)}(z)$ and $R^{(k,m)}(z)$ are 
well defined as operators acting on
$V^{(m)} \otimes V^{(k)}$ and $V^{(k)} \otimes V^{(m)}$ 
respectively.
\end{proposition}
\vspace{5mm}

Operators $R^{(m,k)}(z)$ and $R^{(k,m)}(z)$ satisfy the intertwining property
\eqref{R-matrix} i.e. they are R-matrices on $V^{(k)}\otimes V^{(m)}$ and $V^{(m)} \otimes V^{(k)}$. 
\begin{proposition}
For all $x \in U_q(\frak g)$,
\begin{eqnarray}
R^{(m,k)}(z_1/z_2)
(\pi_{z_1}^{(m)} \otimes \pi_{z_2}^{(k)}) \bigtriangleup(x) 
&=& (\pi_{z_1}^{(m)} \otimes \pi_{z_2}^{(k)}) \bigtriangleup'(x) 
R^{(m,k)}(z_1/z_2),\label{fusion3} \\
R^{(k,m)}(z_1/z_2)
(\pi_{z_1}^{(k)} \otimes \pi_{z_2}^{(m)}) \bigtriangleup(x) 
&=& (\pi_{z_1}^{(k)} \otimes \pi_{z_2}^{(m)}) \bigtriangleup'(x) 
R^{(k,m)}(z_1/z_2),\label{fusion4} 
\end{eqnarray}
\end{proposition}
\noindent
Proof. \ Here we only prove \eqref{fusion3}.
We denote
\[
 \bigtriangleup(x) = \sum x_{(1)} \otimes x_{(2)},
\]
\[
 (\bigtriangleup \otimes \operatorname{id})(x) 
= (\operatorname{id} \otimes \bigtriangleup)(x) =
\sum x_{(1)} \otimes x_{(2)} \otimes x_{(3)}.
\]
For all $x \in U_q(\frak g)$,
\begin{equation*} 
\begin{split}  
 &R^{(m,k)}(z_1/z_2)
(\pi_{z_1}^{(m)} \otimes \pi_{z_2}^{(k)}) \bigtriangleup(x) \\  
 &\;\; =R^{(m,k)}(z_1/z_2)
      \{ \sum \pi_{z_1}^{(m)}(x_{(1)}) \otimes \pi_{z_2}^{(k)}(x_{(2)}) \}\\  
 &\;\; ={\bar{R}}^{(m,n')}((-q)^{-n+k+1}z_1/z_2)_{13}
   {\bar{R}}^{(m,n)}((-q)^{n-k-1}z_1/z_2)_{12} \\  
 &\;\;  \quad \times \{ \sum \pi_{z_1}^{(m)}( x_{(1)}) \otimes 
           \pi_{(-q)^{-n+k+1}z_2}^{(n)}( x_{(2)} )\otimes 
           \pi_{(-q)^{n-k-1}z_2}^{(n')}( x_{(3)}) \}\\  
 &\;\; ={\bar{R}}^{(m,n')}((-q)^{-n+k+1}z_1/z_2)_{13}
      \{\sum \pi_{z_1}^{(m)}( x_{(2)}) \otimes 
           \pi_{(-q)^{-n+k+1}z_2}^{(n)}( x_{(1)} )\otimes 
           \pi_{(-q)^{n-k-1}z_2}^{(n')}( x_{(3)})\} \\  
 &\;\;  \quad \times {\bar{R}}^{(m,n)}((-q)^{n-k-1}z_1/z_2)_{12}\\ 
 &\;\; =\{\sum \pi_{z_1}^{(m)}( x_{(3)}) \otimes 
           \pi_{(-q)^{-n+k+1}z_2}^{(n)}( x_{(2)} )\otimes 
           \pi_{(-q)^{n-k-1}z_2}^{(n')}( x_{(1)}) \} \\
 &\;\; \quad \times {\bar{R}}^{(m,n')}((-q)^{-n+k+1}z_1/z_2)_{13}
       {\bar{R}}^{(m,n)}((-q)^{n-k-1}z_1/z_2)_{12} \\  
 &\;\; =\{\sum \pi_{z_1}^{(m)}(x_{(2)}) \otimes \pi_{z_2}^{(k)}(x_{(1)}) \}
       R^{(m,k)}(z_1/z_2) \\  
 &\;\; =(\pi_{z_1}^{(m)}\otimes \pi_{z_2}^{(k)}) \bigtriangleup'(x) 
       R^{(m,k)}(z_1/z_2). 
\end{split}
\end{equation*}
\hfill Q.E.D
\vspace{5mm}

For our aim of finding commutation relations of vertex operators,
it is enough to consider R-matrices in the cases listed below (see Section 4.3):
\begin{flushleft}
  \begin{tabular}{p{20mm}p{30mm}p{20mm}p{20mm}}
     I & $V^{(k)} \otimes V^{(n)}$ & $n$:even & $k$:even \\ 
     II & $V^{(k)} \otimes V^{(n)}$ & $n$:odd  & $k$:even \\
     III & $V^{(k)} \otimes V^{(n-1)}$ & $n$:even & $k$:odd  \\
     IV & $V^{(k)} \otimes V^{(n-1)}$ & $n$:odd  & $k$:odd  \\
     V & $V^{(n)} \otimes V^{(k)}$ & $n$:even & $k$:even \\
     VI & $V^{(n)} \otimes V^{(k)}$ & $n$:odd  & $k$:even \\
     VII & $V^{(n-1)} \otimes V^{(k)}$ & $n$:even & $k$:odd  \\
     VIII & $V^{(n-1)} \otimes V^{(k)}$ & $n$:odd  & $k$:odd
  \end{tabular}  
\end{flushleft}
where $1 \leq k \leq n-2$ in all these cases.

Let $u^{(k)}$ be the highest weight vector of the fundamental representation $V^{(k)}$. 
We normalize R-matrices such that the eigenvalues of elements $u^{(k)} \otimes u^{(m)}$
(resp. $u^{(m)} \otimes u^{(k)}$) for $\bar{R}^{(k,m)}_n(z)$ (resp. $\bar{R}^{(m,k)}_n(z)$)
are equal to one, then we have

\begin{center}
  \begin{tabular}{cc}
I and II \quad & $ \qquad \bar{R}^{(k,n)}_n(z) = {\displaystyle \prod_{i=1}^{[ \frac{n-k}{2} ]}}
 a(q^{4i-4}(-q)^{-n+k+1}z)^{-1} R^{(k,n)}_n(z),$ \\
\; & \; \\
III and IV \quad & $ \qquad \bar{R}^{(k,n-1)}_n(z) ={\displaystyle \prod_{i=1}^{[ \frac{n-k-1}{2} ]}}
 a(q^{4i-2}(-q)^{-n+k+1}z)^{-1} R^{(k,n-1)}_n(z),$ \\ 
\; & \; \\
V and VII\quad & $ \qquad \bar{R}^{(m,k)}_n(z) ={\displaystyle \prod_{i=1}^{[ \frac{n-k}{2} ]}}
 a(q^{4i-4}(-q)^{-n+k+1}z)^{-1} R^{(m,k)}_n(z),$ \\
\; & \; \\
VI and VIII\quad & $ \qquad \bar{R}^{(m,k)}_n(z) ={\displaystyle \prod_{i=1}^{[ \frac{n-k-1}{2} ]}}
 a(q^{4i-2}(-q)^{-n+k+1}z)^{-1} R^{(m,k)}_n(z),$ \\
\; & \; \\
   \end{tabular}
\end{center}
where $m = \varphi^{(n)}(k)$ and $a(z) = q(1-z)/(1-q^2z)$.

We denote the scalar functions $\alpha^{V^{(k)},V^{(m)}}(z)$,
$\alpha^{V^{(m)},V^{(k)}}(z)$, $\beta^{V^{(k)},V^{(m)}}(z)$ and $\beta^{V^{(m)},V^{(k)}}(z)$
by $\alpha^{(k,m)}_n(z)$, $\alpha^{(m,k)}_n(z)$, 
$\beta^{(k,m)}_n(z)$ and $\beta^{(m,k)}_n(z)$.
In a similar way to Section 2.6, we can show
\begin{equation*}
 \begin{split}
   &\alpha^{(k,m)}_n(z) = \alpha^{(m,k)}_n(z)
     = \frac{(1+(-1)^{n-k}q^{k}{\xi}^{-1/2}z)(1+(-1)^{n-k}q^{-k}{\xi}^{-3/2}z)}{
        (1+(-1)^{n-k}q^{-k}{\xi}^{-1/2}z)(1+(-1)^{n-k}q^{k}{\xi}^{-3/2}z)},   \\
   &\;\\
   &\beta^{(k,m)}_n(z) = \beta^{(m,k)}_n(z)
     = q^{-\frac{k}{2}}
     \frac{(-(-1)^{n-k}q^{k}{\xi}^{1/2}z;{\xi}^2)_{\infty}(-(-1)^{n-k}q^{-k}{\xi}^{3/2}z;{\xi}^2)_{\infty}}{
      (-(-1)^{n-k}q^{-k}{\xi}^{1/2}z;{\xi}^2)_{\infty}(-(-1)^{n-k}q^{k}{\xi}^{3/2}z;{\xi}^2)_{\infty}}.
 \end{split}
\end{equation*}

\section{Vertex operators and two point functions}

In this section we define vertex operators and calculate two point functions.

\subsection{Vertex operators}

Vertex operators are the following homomorphisms of $U_q(\frak g)$-modules:
\begin{eqnarray*}
  \text{(type I)}\;\; \tilde{\Phi}_{\lambda}^{\mu \; V^{(k)}}(z) \; 
    : \; V(\lambda) \longrightarrow \hat{V}(\mu) \otimes V_z^{(k)},\\
  \text{(type II)}\;\; \tilde{\Phi}_{\lambda}^{V^{(k)} \; \mu}(z) \; 
    : \; V(\lambda) \longrightarrow V_z^{(k)} \otimes \hat{V}(\mu),
\end{eqnarray*}
where $\hat{V}(\mu)$ denotes the completion of $V(\mu)$ in the formal topology
with respect to the weight decomposition of $V(\mu)$. (cf. \cite{DJO})
The spin representations have level one prefect crystal (cf. \cite{KKMMNN}),
then our model can be constructed by using
the vertex operators with level one dominant integral weight $\lambda$ and $\mu$.
Here we remark that the level one dominant integral weights are
only $\Lambda_0,\;\Lambda_1,\;\Lambda_{n-1}$ and $\Lambda_n$ in $D_n^{(1)}$ case.
By using conformal weight $\Delta_\lambda = (\lambda|\lambda+2\rho)/2(h^\vee +1)$ 
(for the level one dominant integral weights we have
$\Delta_{\Lambda_0}=0$, $\Delta_{\Lambda_1}=1/2$, $\Delta_{\Lambda_{n-1}}=n/8$, $\Delta_{\Lambda_n}=n/8$),
we multiply $z^{\Delta_\mu-\Delta_\lambda}$ to vertex operators
as follows:
\begin{eqnarray*}
  \Phi_{\lambda}^{\mu \; V^{(k)}}(z)
    =z^{\Delta_\mu-\Delta_\lambda}\tilde{\Phi}_{\lambda}^{\mu \; V^{(k)}}(z),\\
  \Phi_{\lambda}^{V^{(k)} \; \mu}(z)
    =z^{\Delta_\mu-\Delta_\lambda}\tilde{\Phi}_{\lambda}^{V^{(k)} \; \mu}(z).
\end{eqnarray*} 

When we express the image of the highest weight vector $|\lambda \rangle \in V(\lambda)$ 
under $\tilde{\Phi}_{\lambda}^{\mu \; V^{(k)}}(z)$ as a linear combination of the weight vectors of $V(\mu)$,
we call the coefficient of $|\mu \rangle$ the leading term of $\tilde{\Phi}_{\lambda}^{\mu \; V^{(k)}}(z)$ i.e.
\[
\tilde{\Phi}_{\lambda}^{\mu \; V^{(k)}}(z) |\lambda \rangle =  |\mu \rangle \otimes v + \cdots,
\]
where $v \in V^{(k)}$ is the leading term. Here we put
\begin{equation*}
  (V^{(k)})_{\lambda}^{\mu} = 
    \{v \in V^{(k)} | \lambda \equiv \mu + \operatorname{wt}v 
        \pmod{\delta}, \;  {e_i}^{\langle h_i,\mu \rangle + 1}v = 0, \; \text{for all}\; i \}.
\end{equation*}
From the following proposition, we can know when non-trivial vertex operator exists.
\begin{proposition}\cite{DJO}
The mapping to send a vertex operator to its leading term gives an isomorphism
of vector space:
\begin{equation*}
  \{ \tilde{\Phi}_{\lambda}^{\mu \; V^{(k)}}(z) \; 
    : \; V(\lambda) \longrightarrow \hat{V}(\mu) \otimes V_z^{(k)} \}
    \simeq (V^{(k)})_{\lambda}^{\mu}.
\end{equation*}
\end{proposition}
Then for a given vector $v \in (V^{(k)})_{\lambda}^{\mu}$,
a vertex operator which satisfies
\[
\tilde{\Phi}_{\lambda}^{\mu \; V^{(k)}}(z) |\lambda \rangle =  |\mu \rangle \otimes v + \cdots,
\]
uniquely exists. We will normalize vertex operators by specifying the leading term $v$.

In the case of $U_q(D_n^{(1)})$, 
\begin{equation}
  \operatorname{dim} (V^{(k)})_{\lambda}^{\mu} = 0 \; \text{or} \; 1.
  \label{uniqueness}
\end{equation}
Therefore if there exists  non-trivial vertex operators then it is unique
up to multiple constant.
Only in the following cases, non-trivial vertex operator 
$\tilde{\Phi}_{\lambda}^{\mu \; V^{(k)}}(z)$ exists. 
The existence and uniqueness conditions for type II vertex operators are exactly same (cf. \cite{DaiO}):

\newpage
\begin{table}[h]
\begin{center}
\label{tab:tab1}
  \begin{tabular}{| c | c | c |}
     \hline
     $\lambda$  & $\mu$ & $k$ \\
     \hline
     $\Lambda_{n}$   & $\Lambda_{0}$   & $n$ \\ 
     $\Lambda_{n}$   & $\Lambda_{1}$   & $n-1$ \\
     $\Lambda_{n-1}$ & $\Lambda_{0}$   & $n-1$ \\
     $\Lambda_{n-1}$ & $\Lambda_{1}$   & $n$ \\
     $\Lambda_{0}$   & $\Lambda_{n}$   & $n-i$ \\
     $\Lambda_{0}$   & $\Lambda_{n-1}$ & $n+i-1$ \\
     $\Lambda_{1}$   & $\Lambda_{n}$   & $n+i-1$ \\
     $\Lambda_{1}$   & $\Lambda_{n-1}$ & $n-i$ \\
     \hline
  \end{tabular} 
  \qquad
  \begin{tabular}{| c | c | c |}
     \hline
     $\lambda$  & $\mu$ & $k$ \\
     \hline
     $\Lambda_{0}$   & $\Lambda_{0}$   & $1 \leq k \leq n-2$ : even \\ 
     $\Lambda_{1}$   & $\Lambda_{0}$   & $1 \leq k \leq n-2$ : odd \\
     $\Lambda_{1}$   & $\Lambda_{1}$   & $1 \leq k \leq n-2$ : even \\
     $\Lambda_{0}$   & $\Lambda_{1}$   & $1 \leq k \leq n-2$ : odd \\
     $\Lambda_{n}$   & $\Lambda_{n}$   & $1 \leq k \leq n-2$ : even \\
     $\Lambda_{n-1}$ & $\Lambda_{n}$   & $1 \leq k \leq n-2$ : odd \\
     $\Lambda_{n-1}$ & $\Lambda_{n-1}$ & $1 \leq k \leq n-2$ : even \\
     $\Lambda_{n}$   & $\Lambda_{n-1}$ & $1 \leq k \leq n-2$ : odd \\
     \hline
  \end{tabular}  
\end{center}
where $i=0$ if $n$ is even and $i=1$ if $n$ is odd.  
\end{table}

After the normalization of three vertex operators
$\tilde{\Phi}^{\Lambda_{0}\;V^{(n)}}_{\Lambda_{n}}(z)$, 
$\tilde{\Phi}^{\Lambda_{0}\;V^{(2k)}}_{\Lambda_{0}}(z)$ and $\tilde{\Phi}^{\Lambda_{0}\;V^{(2k-1)}}_{\Lambda_{1}}(z)$,
the normalization of the others can be determined
by using Dynkin diagram automorphisms in the same way as \cite{DaiO}.
Here we can choose arbitrary normalization of the above three vertex operators,
and so we fix the normalization as follows:
\begin{eqnarray*}
  \tilde{\Phi}^{\Lambda_{0}\;V^{(n)}}_{\Lambda_{n}}(z)|\Lambda_n \rangle
    &=& |\Lambda_0 \rangle \otimes v_1 + \cdots,\\
  \tilde{\Phi}^{\Lambda_{0}\;V^{(k)}}_{\Lambda_{0}}(z)|\Lambda_0 \rangle
    &=& |\Lambda_0 \rangle \otimes v_2 + \cdots\quad(k:\text{even}),\\
  \tilde{\Phi}^{\Lambda_{0}\;V^{(k)}}_{\Lambda_{1}}(z)|\Lambda_1 \rangle
    &=& |\Lambda_0 \rangle \otimes v_3 + \cdots\quad(k:\text{odd}).
\end{eqnarray*}
The leading terms $v_1$, $v_2$ and $v_3$ are given by
\begin{equation}
  \begin{split}
  v_1 &= v_{(+,+, \cdots,+)}, \\
  v_2 
    &= \sum_{m=0}^{ \frac{k}{2} } c_m^{(n,k)} u_n^{(n-k+2m)}, \\
  v_3 
    &= \sum_{m=0}^{ \frac{k-1}{2} } c_m^{(n,k)} v_{++} \otimes u_{n-1}^{(n-k+2m)},
  \end{split}
  \label{leading term}
\end{equation}
where 
\[
  c_m^{(n,k)}
    = \prod_{j=0}^{m-1}
      \frac{q^{-2(j+1)}(1+q^{2m-2j})}{(1+q^{2n-2k+4j})(1+q^{2n-2k+4j+2})(1-q^{2n-2k+2j})},
\]
and $u_n^{(k)}$ is inductively defined by

\vspace{5mm}
\noindent for $ n \geq 2 $,
\begin{eqnarray*}
  & u_n^{(k)}
      =-\dfrac{(-q)^{-n+k+2}}{(1+q^{2k})(1+q^{2k+2})}
        \{v_{+-} \otimes u_{n-1}^{(k+1)}
           -(-q)^{-k-1}v_{-+}\otimes \sigma_{n-1} \tilde{u}_{n-1}^{(k+1)} \}& \\
  &\quad \qquad +(-q)^{-n+k}\{v_{+-} \otimes u_{n-1}^{(k-1)}
            +(-q)^{k-1}v_{-+}\otimes \sigma_{n-1} \tilde{u}_{n-1}^{(k-1)} \}&
        \hspace{-5mm} (1 \leq k \leq n),\\
  & \tilde{u}_n^{(k)}
      =-\dfrac{(-q)^{-n+k+2}}{(1+q^{2k})(1+q^{2k+2})}
        \{v_{+-} \otimes \tilde{u}_{n-1}^{(k+1)}
           -(-q)^{-k-1}v_{-+}\otimes \sigma_{n-1} u_{n-1}^{(k+1)} \}& \\
  & \quad \qquad +(-q)^{-n+k}\{v_{+-} \otimes \tilde{u}_{n-1}^{(k-1)}
            +(-q)^{k-1}v_{-+}\otimes \sigma_{n-1} u_{n-1}^{(k-1)}\}&
        \hspace{-5mm} (1 \leq k \leq n), \\ 
  & \hspace{-25mm} u_n^{(0)} 
      =-\dfrac{(-q)^{-n+2}}{(1+q^2)}
      \{v_{+-} \otimes u_{n-1}^{(1)}
       -(-q)^{-1}v_{-+}\otimes \sigma_{n-1} \tilde{u}_{n-1}^{(1)} \} &  (k = 0),\\
  & \hspace{-114mm} \tilde{u}_n^{(0)}=0 & (k = 0),\\
  & \hspace{-101mm} u_n^{(k)}=\tilde{u}_n^{(k)} = 0 & \hspace{-10mm} (k > n \; \text{or} \; k < 0),
\end{eqnarray*}

\noindent for $n = 1$,
\begin{eqnarray*} 
  u_1^{(k)}=\tilde{u}_1^{(k)} = v_{+-} \;(k = 1), \quad 
u_1^{(k)}=\tilde{u}_1^{(k)} =0 \;(k \not= 1).
\end{eqnarray*} 

The normalization of type II vertex operators are  given by
\begin{eqnarray*}
  \tilde{\Phi}^{V^{(n)}\;\Lambda_{0}}_{\Lambda_{n}}(z)|\Lambda_n \rangle
    &=&  v_1 \otimes |\Lambda_0 \rangle + \cdots,\\
  \tilde{\Phi}^{V^{(k)}\;\Lambda_{0}}_{\Lambda_{0}}(z)|\Lambda_0 \rangle
    &=&  v_2 \otimes |\Lambda_0 \rangle + \cdots\quad(k:\text{even}),\\
  \tilde{\Phi}^{V^{(k)}\;\Lambda_{0}}_{\Lambda_{1}}(z)|\Lambda_1 \rangle
    &=&  v_3 \otimes |\Lambda_0 \rangle + \cdots\quad(k:\text{odd}),
\end{eqnarray*}
where the leading term $v_1$, $v_2$ and $v_3$ are the same ones in \eqref{leading term}.
The normalization of the other cases is given by Dynkin diagram automorphisms.

\subsection{Vertex operators for dual modules}
By using the isomorphisms $C^{(k)}_{\pm} \; (k = 1, 2, \cdots, n)$,
we can construct vertex operators for dual modules.
These isomorphisms are defined in \eqref{dual for fusion} for $k = 1, 2, \cdots, n-2$ and 
$C^{(n-1)}_{\pm}=C^{(-)}_{\pm}$ and $C^{(n)}_{\pm}=C^{(+)}_{\pm}$ (see \eqref{dual for spin}).
Here we denote by $k'$ the number such that $V^{(k')}_{z \xi^{\mp 1}}$ is isomorphic to 
$(V^{(k)}_z)^{*a^{\pm1}}$ as $U_q(\frak g)$-module, i.e.
\begin{eqnarray*}
  k' = 
  \begin{cases}
   k \quad   &(1 \leq k \leq n-2), \\
   n-1 \quad &(n:\text{even}, \; k = n-1), \\
   n \quad &(n:\text{even}, \; k = n), \\
   n \quad &(n:\text{odd}, \; k = n-1), \\
   n-1 \quad &(n:\text{odd}, \; k = n).
  \end{cases}
\end{eqnarray*}

Thanks to Proposition 4.4 and \eqref{uniqueness},
the fact that $V^{(k')}_{z \xi^{\mp 1}}$ are isomorphic to $(V^{(k)}_z)^{*a^{\pm1}}$
leads that if there exists the vertex operator $\Phi^{\mu \; (V^{(k)})^{*a^{\pm1}}}_{\lambda}(z)$
then it is unique up to multiple scalar factor.
Therefore we have
\begin{equation*}
  \begin{split}
  \tilde{\Phi}^{\mu \; (V^{(k)})^{*a^{\pm1}}}_{\lambda}(z)
    = (const.) (\operatorname{id} \otimes C^{(k')}_{\pm}) 
       \tilde{\Phi}^{\mu \; V^{(k')}}_{\lambda}(z \xi^{\mp 1}),\\
  \tilde{\Phi}^{(V^{(k)})^{*a^{\pm1}} \; \mu}_{\lambda}(z)
    = (const.) (C^{(k')}_{\pm} \otimes \operatorname{id}) 
       \tilde{\Phi}^{V^{(k')} \; \mu}_{\lambda}(z \xi^{\mp 1}).
  \end{split}
\end{equation*} 

In order to construct the transfer matrix and the creation operator, 
it is enough to consider type I vertex operators
$\Phi^{\mu \; (V^{(k)})^{*a^{\pm1}}}_{\lambda}(z)$ for $k=n$ and type II vertex operators
$\Phi^{(V^{(k)})^{*a^{-1}} \; \mu}_{\lambda}(z)$ for $k=1,\ldots,n$. (see Section 6)

Let $\{(v_{(\varepsilon_1, ... ,\varepsilon_n)})^* \}$
be dual basis of $\{v_{(\varepsilon_1, ... ,\varepsilon_n)} \}$.
We give the normalization of the vertex operator 
$\tilde{\Phi}^{\Lambda_0 \; (V^{(n)})^{*a^{\pm1}}}_{\Lambda_n}(z)$ by
\begin{equation*}
  \tilde{\Phi}^{\Lambda_0 \; (V^{(n)})^{*a^{\pm1}}}_{\Lambda_n}(z) |\Lambda_n \rangle
     = |\Lambda_0 \rangle \otimes (v_{(-,-,\cdots,-)})^* + \cdots,
\end{equation*}
and the other vertex operators are normalized by an appropriate Dynkin Diagram automorphism
$\hat{\sigma}$ as follows:
\begin{equation*}
  \tilde{\Phi}^{\hat{\sigma}(\Lambda_0) \; (\hat{\sigma}(V^{(n)}))^{*a^{\pm1}}}_{\hat{\sigma}(\Lambda_n)}(z) 
     |\hat{\sigma}(\Lambda_n) \rangle
     = |\hat{\sigma}(\Lambda_0) \rangle \otimes (\hat{\sigma} v_{(-,-,\cdots,-)})^* + \cdots.
\end{equation*}

We normalize vertex operators of type II such that
\begin{equation*}
  \tilde{\Phi}^{(V^{(k)})^{*a^{-1}} \; \mu}_{\lambda}(z)
    = s (\operatorname{id} \otimes C^{(k')}_{-}) 
       \tilde{\Phi}^{V^{(k')} \; \mu}_{\lambda}(z \xi),
\end{equation*}
where $s$ is a constant defined by
\begin{eqnarray*}
    s =
     \begin{cases}
     (-1)^{[ \frac{n}{2} ] + k(n-k)} \quad 
    &( \lambda = \Lambda_0, \; \Lambda_1 \; \text{and} \; k = n-1,n ),\\
     (-1)^{k(n-k)}
    &( \lambda = \Lambda_0, \; \Lambda_1 \; \text{and} \; k \leq n-2 ),\\
     1 
    &( \text{others}).
    \end{cases}
\end{eqnarray*}
We remark that the constant $s$ is defined 
as above so that the scalar factor $\tau^{(k)}(z)$ in the commutation relation \eqref{commutation relation 3}
does not depend on the fundamental weights $\lambda$ and $\mu$.

\subsection{Two point functions}
We calculate the following vacuum expectation values of 
two vertex operators (two point functions).
\begin{equation}
  \begin{matrix}
  \text{(type I--I)} \hspace{20mm} & {\displaystyle \langle \Phi^{\nu \;W_2}_{\mu}(z_2) 
               \Phi^{\mu \;V_1}_{\lambda}(z_1) \rangle }& \hspace{30mm} \\
  \text{(type II--I)} \hspace{20mm} & {\displaystyle \langle \Phi^{W_2 \; \nu}_{\mu}(z_2)
                 \Phi^{\mu \;V_1}_{\lambda}(z_1) \rangle }& \hspace{30mm} \\
  \text{(type I--II)} \hspace{20mm} & {\displaystyle \langle \Phi^{\nu \;W_2}_{\mu'}(z_2)
                \Phi^{V_1 \; \mu'}_{\lambda}(z_1) \rangle }& \hspace{30mm}
  \end{matrix}
\end{equation}
In the case of type I-I, consider the composition
\begin{equation*}
  V(\lambda) \overset{\Phi^{\mu \;V}_{\lambda}}{\longrightarrow}
  \hat{V}(\mu) \otimes V_{z_1}
  \overset{\Phi^{\nu \;W}_{\mu} \otimes \operatorname{id}}{\longrightarrow}
  \hat{V}(\nu) \otimes W_{z_2} \otimes V_{z_1}
  \overset{\operatorname{id} \otimes P}{\longrightarrow}
  \hat{V}(\nu) \otimes V_{z_1} \otimes W_{z_2}.
\end{equation*}
Expressing the image of the highest weight vector of $V(\lambda)$
under the above composition as a linear combination of weight vectors
with coefficients in $W_{z_1} \otimes V_{z_2}$, the coefficient of the highest weight
vector of $V(\nu)$ is called vacuum expectation value.
Here the subscript of the space $V^{(k)}$ and $V^{(n)}$ in two point function means that
the space indexed by 1 (resp.2) always come in the first (resp. second) component of tensor product. 
Therefore in the case of type II-I and I-II we do not need the last transposition. (cf. \cite{IIJMNT})

Let $R^{V,W}_+(z)$ be the image of the modified R-matrix in \eqref{modified R-matrix}.
\begin{proposition}
\cite{IIJMNT} Let $\Psi(z_1,z_2)$ be a two point function of type I-I ,II-I or I-II. Then 
$\Psi(z_1,z_2)$ satisfies the following difference equation (the q-KZ equation):
\begin{equation}
  \begin{split}
   &\Psi(pz_1,z_2) = A(z_1/z_2)\Psi(z_1,z_2), \\
   &\Psi(pz_1,pz_2) = (q^{-\phi} \otimes q^{-\phi})\Psi(z_1,z_2),
  \end{split}
  \label{q-KZ equation}
\end{equation}
where $p = q^{2(h^{\vee}+1)}$, $\phi = \bar{\lambda} + \bar{\nu} + 2\bar{\rho}$ and
\begin{eqnarray}
  A(z) &=& R^{V,W}_+(pz)(q^{-\phi} \otimes 1) \hspace{39mm} \text{for type I-I},
     \label{q-KZ I-I}\\
       &=& (q^{-\bar{\nu}} \otimes 1)R^{V,W}_+(pq^{-1}z)(q^{-\phi+\bar{\nu}} \otimes 1) 
           \hspace{10mm} \text{for type II-I}, 
     \label{q-KZ II-I}\\
       &=&(q^{-\phi+\bar{\nu}} \otimes 1)R^{V,W}_+(qz)(q^{-\bar{\nu}} \otimes 1) 
           \hspace{16mm} \text{for type I-II}.
     \label{q-KZ I-II}
\end{eqnarray}
\end{proposition}
\noindent Then we can determine two point functions as solutions of the q-KZ equation.

For our aim to determine commutation relations of vertex operators, 
we need to calculate the following four types of two point functions:
\begin{equation}
  \begin{matrix}
  \text{(i)} \hspace{25mm} &{\displaystyle \langle \Phi^{\nu \;(V^{(k)})_2}_{\mu}(z_2)
               \Phi^{\mu \;(V^{(n)})_1}_{\lambda}(z_1) \rangle }& \hspace{25mm}\\
  \text{(ii)} \hspace{25mm} &{\displaystyle \langle \Phi^{\nu \;(V^{(n)})_2}_{\mu'}(z_2)
                \Phi^{\mu' \;(V^{(k)})_1}_{\lambda}(z_1) \rangle }& \hspace{25mm}\\
  \text{(iii)} \hspace{25mm} &{\displaystyle \langle \Phi^{\nu \;(V^{(n)})_2}_{\mu'}(z_2)
                \Phi^{(V^{(k)})_1 \; \mu'}_{\lambda}(z_1) \rangle }& \hspace{25mm}\\
  \text{(iv)} \hspace{25mm} &{\displaystyle \langle \Phi^{(V^{(k)})_1 \; \nu}_{\mu}(z_1)
                 \Phi^{\mu \;(V^{(n)})_2}_{\lambda}(z_2) \rangle }& \hspace{25mm}\\
  \end{matrix}
  \label{list of two point function}
\end{equation} 
where $k = 1,2,\cdots,n$, and all combinations of weights $(\nu, \mu', \mu, \lambda)$
that non-trivial vertex operator exists are given by\\
For $n$: even
\begin{center}
  {\bf Table 1}

  \begin{tabular}{| c | c | c | c | c |}
     \hline
     $\nu$ & $\mu'$ & $\mu$ & $\lambda$ & $k$ \\
     \hline
     $\Lambda_n$ & $\Lambda_0$ & $\Lambda_0$ & $\Lambda_n$ & $n$ \\
     $\Lambda_{n-1}$ & $\Lambda_1$ & $\Lambda_1$ & $\Lambda_{n-1}$ & $n$ \\
     $\Lambda_0$ & $\Lambda_n$ & $\Lambda_n$ & $\Lambda_0$ & $n$ \\
     $\Lambda_1$ & $\Lambda_{n-1}$ & $\Lambda_{n-1}$ & $\Lambda_1$ & $n$\\                   
     \hline
  \end{tabular} 
\end{center}
\begin{center}
  \begin{tabular}{| c | c | c | c | c |}
     \hline
     $\nu$ & $\mu'$ & $\mu$ & $\lambda$ & $k$ \\
     \hline
     $\Lambda_n$ & $\Lambda_0$ & $\Lambda_1$ & $\Lambda_{n-1}$ & $n-1$\\
     $\Lambda_{n-1}$ & $\Lambda_1$ & $\Lambda_0$ & $\Lambda_n$ & $n-1$\\
     $\Lambda_0$ & $\Lambda_n$ & $\Lambda_{n-1}$ & $\Lambda_1$ & $n-1$\\
     $\Lambda_1$ & $\Lambda_{n-1}$ & $\Lambda_n$ & $\Lambda_0$ & $n-1$\\
     \hline
 \end{tabular} 
\end{center}
\begin{center}
  \begin{tabular}{| c | c | c | c | c |}
     \hline
     $\nu$ & $\mu'$ & $\mu$ & $\lambda$ & $k$ \\
     \hline
     $\Lambda_n$ & $\Lambda_0$ & $\Lambda_n$ & $\Lambda_0$ & $1 \leq k \leq n-2$ :even \\
     $\Lambda_{n-1}$ & $\Lambda_1$ & $\Lambda_{n-1}$ & $\Lambda_1$ & $1 \leq k \leq n-2$ :even \\
     $\Lambda_0$ & $\Lambda_n$ & $\Lambda_0$ & $\Lambda_n$ & $1 \leq k \leq n-2$ :even \\
     $\Lambda_1$ & $\Lambda_{n-1}$ & $\Lambda_1$ & $\Lambda_{n-1}$ & $1 \leq k \leq n-2$ :even\\
     \hline
  \end{tabular} 
\end{center}
\begin{center}
  \begin{tabular}{| c | c | c | c | c |}
     \hline
     $\nu$ & $\mu'$ & $\mu$ & $\lambda$ & $k$ \\
     \hline
     $\Lambda_n$ & $\Lambda_0$ & $\Lambda_{n-1}$ & $\Lambda_1$ & $1 \leq k \leq n-2$ :odd\\
     $\Lambda_{n-1}$ & $\Lambda_1$ & $\Lambda_n$ & $\Lambda_0$ & $1 \leq k \leq n-2$ :odd \\
     $\Lambda_0$ & $\Lambda_n$ & $\Lambda_1$ & $\Lambda_{n-1}$ & $1 \leq k \leq n-2$ :odd \\
     $\Lambda_1$ & $\Lambda_{n-1}$ & $\Lambda_0$ & $\Lambda_n$ & $1 \leq k \leq n-2$ :odd \\ 
     \hline
\end{tabular} 
\end{center}

For $n$: odd
\begin{center}
  {\bf Table 2}

  \begin{tabular}{| c | c | c | c | c |}
     \hline
     $\nu$ & $\mu'$ & $\mu$ & $\lambda$ & $k$ \\
     \hline
     $\Lambda_n$ & $\Lambda_1$ & $\Lambda_1$ & $\Lambda_{n-1}$ & $n$\\
     $\Lambda_{n-1}$ & $\Lambda_0$ & $\Lambda_0$ & $\Lambda_{n}$ & $n$ \\
     $\Lambda_0$ & $\Lambda_{n}$ & $\Lambda_{n}$ & $\Lambda_1$ & $n$ \\
     $\Lambda_1$ & $\Lambda_{n-1}$ & $\Lambda_{n-1}$ & $\Lambda_0$ & $n$\\ 
     \hline
 \end{tabular} 
\end{center}
\begin{center}
  \begin{tabular}{| c | c | c | c | c |}
     \hline
     $\nu$ & $\mu'$ & $\mu$ & $\lambda$ & $k$ \\
     \hline
     $\Lambda_n$ & $\Lambda_1$ & $\Lambda_0$ & $\Lambda_n$ & $n-1$ \\
     $\Lambda_{n-1}$ & $\Lambda_0$ & $\Lambda_1$ & $\Lambda_{n-1}$ & $n-1$ \\
     $\Lambda_0$ & $\Lambda_n$ & $\Lambda_{n-1}$ & $\Lambda_0$ & $n-1$ \\
     $\Lambda_1$ & $\Lambda_{n-1}$ & $\Lambda_n$ & $\Lambda_1$ & $n-1$ \\
     \hline
 \end{tabular} 
\end{center}
\begin{center}
  \begin{tabular}{| c | c | c | c | c |}
     \hline
     $\nu$ & $\mu'$ & $\mu$ & $\lambda$ & $k$ \\
     \hline
     $\Lambda_n$ & $\Lambda_1$ & $\Lambda_n$ & $\Lambda_1$ & $1 \leq k \leq n-2$ :even \\
     $\Lambda_{n-1}$ & $\Lambda_0$ & $\Lambda_{n-1}$ & $\Lambda_0$ & $1 \leq k \leq n-2$ :even \\
     $\Lambda_0$ & $\Lambda_n$ & $\Lambda_0$ & $\Lambda_n$ & $1 \leq k \leq n-2$ :even \\
     $\Lambda_1$ & $\Lambda_{n-1}$ & $\Lambda_1$ & $\Lambda_{n-1}$ & $1 \leq k \leq n-2$ :even\\ 
     \hline
 \end{tabular} 
\end{center}
\begin{center}
  \begin{tabular}{| c | c | c | c | c |}
     \hline
     $\nu$ & $\mu'$ & $\mu$ & $\lambda$ & $k$ \\
     \hline 
     $\Lambda_n$ & $\Lambda_1$ & $\Lambda_{n-1}$ & $\Lambda_0$ & $1 \leq k \leq n-2$ :odd \\
     $\Lambda_{n-1}$ & $\Lambda_0$ & $\Lambda_n$ & $\Lambda_1$ & $1 \leq k \leq n-2$ :odd  \\
     $\Lambda_0$ & $\Lambda_n$ & $\Lambda_1$ & $\Lambda_{n-1}$ & $1 \leq k \leq n-2$ :odd  \\
     $\Lambda_1$ & $\Lambda_{n-1}$ & $\Lambda_0$ & $\Lambda_n$ & $1 \leq k \leq n-2$ :odd  \\
     \hline
 \end{tabular} 
\end{center}
\vspace{5mm}

The following propositions in \cite{IIJMNT}, \cite{DaiO} are useful to calculate two point functions.
\begin{proposition}\cite{IIJMNT}
Let $\Psi(z_1,z_2)$ be a two point function of type $(i)$ in \eqref{list of two point function}.
For any $i=0,1,\cdots,n$, 
\begin{eqnarray*}
(\pi_{z_1}^{(n)} \otimes \pi_{z_2}^{(k)})
\bigtriangleup'(e_i)^{\langle h_i, \nu \rangle+1}\Psi(z_1, z_2) = 0, \quad
\operatorname{wt} \Psi(z_1, z_2) = \bar{\lambda}-\bar{\nu},
\end{eqnarray*}
hold.
\end{proposition}
\vspace{5mm}
\begin{proposition} \cite{IIJMNT}
Let $\Psi(z_1,z_2)$ be a solution of the q-KZ equation \eqref{q-KZ equation} with $A(z)$ specified by
\eqref{q-KZ I-I}, then 
\begin{equation*}
  \begin{split}
    &(q^{-\bar{\nu}} \otimes 1)\Psi(q^{-1}z_1,z_2), \quad (q^{-\phi+\bar{\nu}} \otimes 1)\Psi(p^{-1}qz_1,z_2),
  \end{split}
\end{equation*}
satisfy the equation with $A(z)$ specified by \eqref{q-KZ II-I} and \eqref{q-KZ I-II} respectively.
\end{proposition}
\vspace{5mm}
\begin{proposition} \cite{DaiO}
If a $V \otimes W$-valued function $w(z)$ satisfies
\begin{equation*}
  \begin{split}
  &(\pi^V_{z_1}\otimes\pi^W_{z_2})\Delta'(e_i)^{\langle h_i,\nu \rangle + 1}w(z_1/z_2) = 0
     \quad (i = 0, 1, \ldots, n),\\
  &\bar{R}^{V,W}(pz)(q^{-\phi} \otimes 1)w(z) = r(z)w(pz),  
  \end{split}
\end{equation*}
for some scalar function $r(z)$, then $\bar{w}(z) = P(q^{-\phi} \otimes 1)w(p^{-1}z^{-1})$
satisfies
\begin{equation*}
  \begin{split}
  &(\pi^W_{z_1}\otimes\pi^V_{z_2})\Delta'(e_i)^{\langle h_i,\nu \rangle + 1}\bar{w}(z_1/z_2) = 0
     \quad (i = 0, 1, \ldots, n),\\
  &\bar{R}^{W,V}(pz)(q^{-\phi} \otimes 1)\bar{w}(z) 
         = q^{\langle \phi,\operatorname{wt}\bar{w} \rangle} r(p^{-2}z^{-1})\bar{w}(pz).  
  \end{split}
\end{equation*}
\end{proposition}
\vspace{5mm}
\begin{proposition} \cite{DaiO}
Let $\lambda$, $\mu$ be level 1 dominant integral weights and $v$ be a weight vector in $V^{(k)}$.
Define a non-negative integer $m(\lambda,\mu;v)$ by the minimal value of $m_0$ satisfying
\begin{eqnarray*}
\lambda - \mu + \sum_{j=0}^{n}m_j\alpha_j \equiv \operatorname{wt}v \pmod{ {\Bbb Z}\delta},\\
m_j \geq 0 \;\;(j = 0,1,\ldots,n). \qquad
\end{eqnarray*}
If a two point function have a form
\[
\langle \Phi^{\nu \;(V^{(\ell)})_2}_{\mu}(z_2)
                          \Phi^{\mu \; (V^{(k)})_1}_{\lambda}(z_1) \rangle
= z_1^{\Delta_{\mu}-\Delta_{\lambda}} z_2^{ \Delta_{\nu}-\Delta_{\mu} }
  \sum_i a_i(z_1/z_2) v_i \otimes v_i',
\]
then 
\begin{eqnarray*}
\langle \Phi^{\hat{\sigma}(\nu) \;\hat{\sigma}(V^{(\ell)})_2}_{\hat{\sigma}(\mu)}(z_2)
            \Phi^{\hat{\sigma}(\mu) \; \hat{\sigma}(V^{(k)})_1}_{\hat{\sigma}(\lambda)}(z_1) \rangle
&=& z_1^{\Delta_{\widehat{\sigma}(\mu)}-\Delta_{\widehat{\sigma}(\lambda)} }
   z_2^{ \Delta_{\widehat{\sigma}(\nu)}-\Delta_{\widehat{\sigma}(\mu)} } \\
& &\times \sum_i a_i(z_1/z_2) (z_1/z_2)^{m_i} \hat{\sigma}(v_i) \otimes \hat{\sigma}(v_i'),
\end{eqnarray*}
where $\hat{\sigma}$ is a Dynkin diagram automorphism and 
\[
m_i = m(\hat{\sigma}(\lambda),\hat{\sigma}(\mu);\hat{\sigma}(v_i)) - m(\lambda,\mu;v_i).
\]
\end{proposition}
\vspace{5mm}

Indeed, by using Proposition 4.4, type (iii) and (iv) in \eqref{list of two point function} 
are determined from type (i) and (ii).
Furthermore we can calculate type (ii) from type (i) by virtue of Proposition 4.5.
Hence we need to know explicit forms of type (i) for each case in Table 1 and Table 2.
Thanks to Dynkin diagram symmetry of two point functions (Proposition 4.6),
we obtain the two point functions in Table 1 and 2 from the ones listed below.
\vspace{5mm}

\noindent For $n$: even
\begin{center}
  \begin{tabular}{| c | c | c | c | c | c |}
     \hline
      \;    & $\nu$ & $\mu'$ & $\mu$ & $\lambda$ & $k$ \\
     \hline
     case 1 & $\Lambda_0$ & $\Lambda_n$ & $\Lambda_n$ & $\Lambda_0$ & $n$ \\
     case 2 & $\Lambda_0$ & $\Lambda_n$ & $\Lambda_{n-1}$ & $\Lambda_1$ & $n-1$\\
     case 3 & $\Lambda_0$ & $\Lambda_n$ & $\Lambda_0$ & $\Lambda_n$ & $1 \leq k \leq n-2$ :even \\   
     case 4 & $\Lambda_0$ & $\Lambda_n$ & $\Lambda_1$ & $\Lambda_{n-1}$ & $1 \leq k \leq n-2$ :odd \\          
     \hline
  \end{tabular} 
\end{center}
\vspace{5mm}
For $n$: odd
\begin{center}
  \begin{tabular}{| c | c | c | c | c | c |}
     \hline
      \;    & $\nu$ & $\mu'$ & $\mu$ & $\lambda$ & $k$ \\
     \hline
     case 5 & $\Lambda_0$ & $\Lambda_{n}$ & $\Lambda_{n}$ & $\Lambda_1$ & $n$ \\
     case 6 & $\Lambda_0$ & $\Lambda_n$ & $\Lambda_{n-1}$ & $\Lambda_0$ & $n-1$ \\
     case 7 & $\Lambda_0$ & $\Lambda_n$ & $\Lambda_0$ & $\Lambda_n$ & $1 \leq k \leq n-2$ :even \\
     case 8 & $\Lambda_0$ & $\Lambda_n$ & $\Lambda_1$ & $\Lambda_{n-1}$ & $1 \leq k \leq n-2$ :odd  \\
     \hline
 \end{tabular} 
\end{center}
Explicit forms of these two point functions are given in the next subsection.

\subsection{Explicit forms of two point functions}

\mbox{}
\vspace{5mm}

\noindent Case 1 

(i) and (ii)
\begin{equation}
    \hspace{-10mm} \langle \Phi^{ {\Lambda}_0 \;(V^{(n)})_2 }_{ {\Lambda}_{n}}(z_2) 
    \Phi^{ {\Lambda}_{n} \; (V^{(n)})_1 }_{ {\Lambda}_0 }(z_1) \rangle 
    = z_1^{\frac{n}{8}} z_2^{-\frac{n}{8}} \psi^{(n,n)}(z_1/z_2) u_n^{(1)} 
\end{equation}

(iii) 
\begin{equation}
  \begin{split}
    & \hspace{-10mm} \langle \Phi^{(V^{(n)})_2 \; {\Lambda}_0 }_{ {\Lambda}_{n}}(z_2)
     \Phi^{ {\Lambda}_{n} \; (V^{(n)})_1 }_{ {\Lambda}_0 }(z_1) \rangle \\
    & \quad = z_1^{\frac{n}{8}} z_2^{-\frac{n}{8}} \psi^{(n,n)}(p^{-1}qz_1/z_2)  
     (-1)^{-\frac{n}{2}}q^{\frac{1}{2}n(n-1)} P u_n^{(1)} 
  \end{split}
\end{equation}

(iv) 
\begin{equation}
  \langle \Phi^{ {\Lambda}_0 \;(V^{(n)})_2 }_{ {\Lambda}_{n}}(z_2)
    \Phi^{ (V^{(n)})_1 \;  {\Lambda}_{n}}_{ {\Lambda}_0 }(z_1) \rangle 
     =  z_1^{-\frac{n}{8}} z_2^{\frac{n}{8}} \psi^{(n,n)}(q^{-1}z_2/z_1) P u_n^{(1)} 
\end{equation}
\vspace{5mm}

\newpage
\noindent Case 2 

(i) 
\begin{equation}
\hspace{-7mm} \langle \Phi^{ {\Lambda}_0 \;(V^{(n-1)})_2 }_{ {\Lambda}_{n-1} }(z_2)
        \Phi^{ {\Lambda}_{n-1} \; (V^{(n)})_1 }_{ {\Lambda}_1 }(z_1) \rangle 
      = z_1^{\frac{n}{8}-\frac{1}{2}} z_2^{-\frac{n}{8}} \psi^{(n,n-1)}(z_1/z_2) u_n^{(2)} 
\end{equation}

(ii)
\begin{equation}
  \begin{split}
    & \hspace{-6mm} \langle \Phi^{ {\Lambda}_0 \;(V^{(n)})_2 }_{ {\Lambda}_{n} }(z_2)
          \Phi^{ {\Lambda}_{n} \; (V^{(n-1)})_1 }_{ {\Lambda}_1 }(z_1) \rangle \\
    & \quad  =  z_1^{\frac{n}{8}-\frac{1}{2}} z_2^{-\frac{n}{8}} \psi^{(n,n-1)}(z_1/z_2) 
        (-1)^{-\frac{n}{2}-1} q^{\frac{1}{2}(2n-1)} (q^{-\phi} \otimes 1) P u_n^{(2)} 
  \end{split}
\end{equation}

(iii)
\begin{equation}
  \begin{split}
    & \hspace{-6mm} \langle \Phi^{ {\Lambda}_0 \;(V^{(n)})_2 }_{ {\Lambda}_{n} }(z_2)
        \Phi^{ (V^{(n-1)})_1 \; {\Lambda}_{n} }_{ {\Lambda}_1 }(z_1) \rangle \\
    & \quad = z_1^{\frac{n}{8}-\frac{1}{2}} z_2^{-\frac{n}{8}}  \psi^{(n,n-1)}(pq^{-1}z_1/z_2)
              (-1)^{-\frac{n}{2}+1}q^{-\frac{1}{2}(n-1)(n-2)}  P u_n^{(2)}
  \end{split}
\end{equation}

(iv)
\begin{equation}
\langle \Phi^{ (V^{(n-1)})_2 \; {\Lambda}_0 }_{ {\Lambda}_{n-1} }(z_2)
        \Phi^{ {\Lambda}_{n-1} \; (V^{(n)})_1 }_{ {\Lambda}_1 }(z_1) \rangle 
  =  z_1^{-\frac{n}{8}} z_2^{\frac{n}{8}-\frac{1}{2}} \psi^{(n,n-1)}(q^{-1}z_2/z_1) P u_n^{(2)}
\end{equation}
\vspace{10mm}

\noindent Case 3 \\

(i)
\begin{equation}
   \hspace{-12mm} \langle \Phi^{{\Lambda}_0 \;(V^{(n)})_2}_{{\Lambda}_{n}}(z_2)
   \Phi^{{\Lambda}_{n} \; (V^{(k)})_1}_{{\Lambda}_n}(z_1) \rangle 
     = z_2^{-\frac{n}{8}} \psi^{(n,k)}(z_1/z_2) w^{(n,k)}(z_1/z_2) 
\end{equation}

(ii)
\begin{equation}
  \begin{split}
   &  \langle \Phi^{{\Lambda}_0 \;(V^{(k)})_2}_{{\Lambda}_{0}}(z_2)
    \Phi^{{\Lambda}_{0} \; (V^{(n)})_1}_{{\Lambda}_n}(z_1) \rangle \\
   & \quad =  z_1^{-\frac{n}{8}} \psi^{(n,k)}(p^{-1}z_2/z_1) 
       q^{\frac{k}{2}(2n-1)} (z_1/z_2)^{\frac{k}{2}} 
          P (q^{-\phi} \otimes 1) w^{(n,k)}(p^{-1}z_2/z_1) 
   \end{split}
\end{equation}

(iii)
\begin{equation}
   \begin{split}
   &  \hspace{-15mm} \langle \Phi^{{\Lambda}_0 \;(V^{(n)})_2}_{{\Lambda}_{n}}(z_2)
     \Phi^{ (V^{(k)})_1 \;{\Lambda}_{n}}_{{\Lambda}_n}(z_1) \rangle \\
   & \quad = z_2^{-\frac{n}{8}} \psi^{(n,k)}(p^{-1}qz_1/z_2) (q^{-\phi} \otimes 1) 
            w^{(n,k)}(p^{-1}qz_1/z_2) 
  \end{split}
\end{equation}

(iv)
\begin{equation}
  \begin{split}
    & \langle \Phi^{(V^{(k)})_1 \;{\Lambda}_0 }_{{\Lambda}_{0}}(z_2)
      \Phi^{{\Lambda}_{0} \; (V^{(n)})_2}_{{\Lambda}_n}(z_1) \rangle \\
    & \quad =   z_2^{-\frac{n}{8}} \psi^{(n,k)}(q^{-1}z_2/z_1)  
       q^{k(n-1)} (z_1/z_2)^{-\frac{k}{2}} (q^{-\phi} \otimes 1) w^{(n,k)}(p^{-1}qz_1/z_2) 
  \end{split}
\end{equation}
\vspace{10mm}

\noindent Case 4 

Two point functions for the case 4 are obtained from the following ones
by using the Dynkin diagram automorphism $\hat{\sigma}_2$.

(i)
\begin{equation}
  \langle \Phi^{{\Lambda}_0 \;(V^{(n-1)})_2}_{{\Lambda}_{n-1}}(z_2)
   \Phi^{{\Lambda}_{n-1} \; (V^{(k)})_1}_{{\Lambda}_n}(z_1) \rangle 
    = z_2^{-\frac{n}{8}} \psi^{(n,k)}(z_1/z_2) w^{(n,k)}(z_1/z_2) 
\end{equation}

(ii)
\begin{equation}
  \begin{split}
    & \hspace{-2mm} \langle \Phi^{{\Lambda}_0 \;(V^{(k)})_2}_{{\Lambda}_{1}}(z_2)
      \Phi^{{\Lambda}_{1} \; (V^{(n-1)})_1}_{{\Lambda}_n}(z_1) \rangle \\
    & \quad = z_1^{\frac{1}{2}-\frac{n}{8}} z_2^{-\frac{1}{2}} \psi^{(n,k)}(p^{-1}z_2/z_1) \\
    & \qquad  \times (-1)^{n-1}q^{\frac{k}{2}(2n-1)} (z_1/z_2)^{\frac{k-1}{2}} P 
              (q^{-\phi} \otimes 1) w^{(n,k)}(p^{-1}z_2/z_1) 
  \end{split}
\end{equation}

(iii)
\begin{equation}
  \begin{split}
    & \hspace{-8mm} \langle \Phi^{{\Lambda}_0 \;(V^{(n-1)})_2}_{{\Lambda}_{n-1}}(z_2)
     \Phi^{ (V^{(k)})_1 \;{\Lambda}_{n-1}}_{{\Lambda}_n}(z_1) \rangle \\
    & \quad  =  z_2^{-\frac{n}{8}} \psi^{(n,k)} q^{\frac{1}{2}} (p^{-1}qz_1/z_2) (q^{-\phi} \otimes 1) 
      w^{(n,k)}(p^{-1}qz_1/z_2) 
  \end{split}
\end{equation}

(iv)
\begin{equation}
  \begin{split}
    & \langle \Phi^{(V^{(k)})_1 \;{\Lambda}_0 }_{{\Lambda}_{1}}(z_2)
      \Phi^{{\Lambda}_{1} \; (V^{(n-1)})_2}_{{\Lambda}_n}(z_1) \rangle \\
    & \quad = z_1^{-\frac{1}{2}} z_2^{\frac{1}{2}-\frac{n}{8}}  \psi^{(n,k)}(q^{-1}z_2/z_1) \\ 
    & \qquad  \times (-1)^{n-1} q^{\frac{1}{2}} q^{k(n-1)} (z_1/z_2)^{-\frac{k-1}{2}}
         (q^{-\phi} \otimes 1) w^{(n,k)}(p^{-1}qz_1/z_2) 
  \end{split}
\end{equation}
\vspace{10mm}

\noindent Case 5

(i) and (ii)
\begin{equation}
  \begin{split}
   \hspace{-7mm} \langle \Phi^{ {\Lambda}_0 \;(V^{(n)})_2 }_{ {\Lambda}_{n}}(z_2) 
    \Phi^{ {\Lambda}_{n} \; (V^{(n)})_1 }_{ {\Lambda}_1 }(z_1) \rangle 
    = z_1^{\frac{n}{8}-\frac{1}{2}} z_2^{-\frac{n}{8}} \psi^{(n,n)}(z_1/z_2) u_n^{(2)} 
  \end{split}
\end{equation}

(iii)
\begin{equation}
  \begin{split}
    & \hspace{-5mm} \langle \Phi^{(V^{(n)})_2 \; {\Lambda}_0 }_{ {\Lambda}_{n}}(z_2)
      \Phi^{ {\Lambda}_{n} \; (V^{(n)})_1 }_{ {\Lambda}_0 }(z_1) \rangle \\
    & \quad = z_1^{\frac{n}{8}-\frac{1}{2}} z_2^{-\frac{n}{8}} \psi^{(n,n)} 
      (-1)^{-\frac{n-1}{2}}q^{\frac{1}{2}(n-1)(n-2)} (p^{-1}qz_1/z_2) P u_n^{(2)} 
  \end{split}
\end{equation}

(iv)
\begin{equation}
 \langle \Phi^{ {\Lambda}_0 \;(V^{(n)})_2 }_{ {\Lambda}_{n}}(z_2)
    \Phi^{ (V^{(n)})_1 \;  {\Lambda}_{n}}_{ {\Lambda}_0 }(z_1) \rangle 
   =  z_1^{-\frac{n}{8}} z_2^{\frac{n}{8}-\frac{1}{2}} \psi^{(n,n)}(q^{-1}z_2/z_1) P u_n^{(2)} 
\end{equation}
\vspace{10mm}
\newpage

\noindent Case 6

(i) 
\begin{equation}
  \hspace{-10mm} \langle \Phi^{ {\Lambda}_0 \;(V^{(n-1)})_2 }_{ {\Lambda}_{n-1} }(z_2)
        \Phi^{ {\Lambda}_{n-1} \; (V^{(n)})_1 }_{ {\Lambda}_0 }(z_1) \rangle 
      = z_1^{\frac{n}{8}} z_2^{-\frac{n}{8}} \psi^{(n,n-1)}(z_1/z_2) u_n^{(1)} 
\end{equation}

(ii)
\begin{equation}
  \begin{split}
      & \hspace{-7mm} \langle \Phi^{ {\Lambda}_0 \;(V^{(n)})_2 }_{ {\Lambda}_{n} }(z_2)
        \Phi^{ {\Lambda}_{n} \; (V^{(n-1)})_1 }_{ {\Lambda}_0 }(z_1) \rangle \\
      & \quad    =  z_1^{\frac{n}{8}} z_2^{-\frac{n}{8}} \psi^{(n,n-1)}(z_1/z_2) 
            (-1)^{-\frac{n-1}{2}} q^{\frac{1}{2}(2n-1)} (q^{-\phi} \otimes 1) P u_n^{(1)} 
  \end{split}
\end{equation}

(iii)
\begin{equation}
  \begin{split}
    & \hspace{-10mm} \langle \Phi^{ {\Lambda}_0 \;(V^{(n)})_2 }_{ {\Lambda}_{n} }(z_2)
        \Phi^{ (V^{(n-1)})_1 \; {\Lambda}_{n} }_{ {\Lambda}_0 }(z_1) \rangle \\
    & \quad =  z_1^{\frac{n}{8}} z_2^{-\frac{n}{8}} \psi^{(n,n-1)}(pq^{-1}z_1/z_2) 
           (-1)^{-\frac{n-1}{2}}q^{-\frac{n}{2}(n-1)} P u_n^{(1)}
  \end{split}
\end{equation}

(iv)
\begin{equation}
   \langle \Phi^{ (V^{(n-1)})_2 \; {\Lambda}_0 }_{ {\Lambda}_{n-1} }(z_2)
       \Phi^{ {\Lambda}_{n-1} \; (V^{(n)})_1 }_{ {\Lambda}_0 }(z_1) \rangle 
     =  z_1^{-\frac{n}{8}} z_2^{\frac{n}{8}} \psi^{(n,n-1)}(q^{-1}z_2/z_1) P u_n^{(1)}
\end{equation}
\vspace{10mm}

\noindent Case 7 

(i)
\begin{equation}
   \hspace{-10mm} \langle \Phi^{{\Lambda}_0 \;(V^{(n)})_2}_{{\Lambda}_{n}}(z_2)
   \Phi^{{\Lambda}_{n} \; (V^{(k)})_1}_{{\Lambda}_n}(z_1) \rangle 
     = z_2^{-\frac{n}{8}} \psi^{(n,k)}(z_1/z_2) w^{(n,k)}(z_1/z_2) 
\end{equation}

(ii)
\begin{equation}
  \begin{split}
   & \langle \Phi^{{\Lambda}_0 \;(V^{(k)})_2}_{{\Lambda}_{0}}(z_2)
    \Phi^{{\Lambda}_{0} \; (V^{(n)})_1}_{{\Lambda}_n}(z_1) \rangle \\
   & \quad  = z_1^{-\frac{n}{8}} \psi^{(n,k)}(p^{-1}z_2/z_1) 
       q^{\frac{k}{2}(2n-1)} (z_1/z_2)^{\frac{k}{2}} 
          P (q^{-\phi} \otimes 1) w^{(n,k)}(p^{-1}z_2/z_1) 
   \end{split}
\end{equation}

(iii)
\begin{equation}
  \begin{split}
   & \hspace{-13mm}\langle \Phi^{{\Lambda}_0 \;(V^{(n)})_2}_{{\Lambda}_{n}}(z_2)
     \Phi^{ (V^{(k)})_1 \;{\Lambda}_{n}}_{{\Lambda}_n}(z_1) \rangle \\
   &  \quad = z_2^{-\frac{n}{8}} \psi^{(n,k)}(p^{-1}qz_1/z_2) (q^{-\phi} \otimes 1) 
              w^{(n,k)}(p^{-1}qz_1/z_2) 
   \end{split}
\end{equation}

(iv)
\begin{equation}
  \begin{split}
   & \langle \Phi^{(V^{(k)})_1 \;{\Lambda}_0 }_{{\Lambda}_{0}}(z_2)
     \Phi^{{\Lambda}_{0} \; (V^{(n)})_2}_{{\Lambda}_n}(z_1) \rangle \\
   & \quad = z_2^{-\frac{n}{8}} \psi^{(n,k)}(q^{-1}z_2/z_1)  
      q^{k(n-1)} (z_1/z_2)^{-\frac{k}{2}} (q^{-\phi} \otimes 1) w^{(n,k)}(p^{-1}qz_1/z_2) 
  \end{split}
\end{equation}
\vspace{10mm}
\newpage

\noindent Case 8 

Two point functions for the case 8 are obtained from the following ones
by using the Dynkin diagram automorphism $\hat{\sigma}_2$.

(i)
\begin{equation}
  \langle \Phi^{{\Lambda}_0 \;(V^{(n-1)})_2}_{{\Lambda}_{n-1}}(z_2)
   \Phi^{{\Lambda}_{n-1} \; (V^{(k)})_1}_{{\Lambda}_n}(z_1) \rangle 
    = z_2^{-\frac{n}{8}} \psi^{(n,k)}(z_1/z_2) w^{(n,k)}(z_1/z_2) 
\end{equation}

(ii)
\begin{equation}
  \begin{split}
   & \hspace{-2mm} \langle \Phi^{{\Lambda}_0 \;(V^{(k)})_2}_{{\Lambda}_{1}}(z_2)
      \Phi^{{\Lambda}_{1} \; (V^{(n-1)})_1}_{{\Lambda}_n}(z_1) \rangle \\
   & \quad = z_1^{\frac{1}{2}-\frac{n}{8}} z_2^{-\frac{1}{2}} \psi^{(n,k)}(p^{-1}z_2/z_1)  \\
   & \qquad  \times  (-1)^{n-1}q^{\frac{k}{2}(2n-1)} (z_1/z_2)^{\frac{k-1}{2}} P 
                     (q^{-\phi} \otimes 1) w^{(n,k)}(p^{-1}z_2/z_1) 
  \end{split}
\end{equation}

(iii)
\begin{equation}
  \begin{split}
   & \hspace{-7mm} \langle \Phi^{{\Lambda}_0 \;(V^{(n-1)})_2}_{{\Lambda}_{n-1}}(z_2)
     \Phi^{ (V^{(k)})_1 \;{\Lambda}_{n-1}}_{{\Lambda}_n}(z_1) \rangle \\
  & \quad  =  z_2^{-\frac{n}{8}} \psi^{(n,k)} q^{\frac{1}{2}} (p^{-1}qz_1/z_2) (q^{-\phi} \otimes 1) 
     w^{(n,k)}(p^{-1}qz_1/z_2) 
  \end{split}
\end{equation}

(iv)
\begin{equation}
  \begin{split}
   & \langle \Phi^{(V^{(k)})_1 \;{\Lambda}_0 }_{{\Lambda}_{1}}(z_2)
     \Phi^{{\Lambda}_{1} \; (V^{(n-1)})_2}_{{\Lambda}_n}(z_1) \rangle \\
   & \quad  = z_1^{-\frac{1}{2}} z_2^{\frac{1}{2}-\frac{n}{8}}  \psi^{(n,k)}(q^{-1}z_2/z_1) \\
    & \qquad  \times  (-1)^{n-1} q^{\frac{1}{2}} q^{k(n-1)} (z_1/z_2)^{-\frac{k-1}{2}}
             (q^{-\phi} \otimes 1) w^{(n,k)}(p^{-1}qz_1/z_2) 
  \end{split}
\end{equation}
\vspace{5mm}
where $\psi^{(n,k)}(z)$ is given by
\begin{equation*}
\psi^{(n,k)}(z) =
  \begin{cases} 
    \displaystyle{\prod^{k-1}_{j=1}} 
    \dfrac{ (-q\xi^{5/2}q^{2j-k+1}(-1)^{n-k}z_1/z_2;\xi^2)_{\infty} }{
          (-q\xi^{3/2}q^{2j-k+1}(-1)^{n-k}z_1/z_2;\xi^2)_{\infty} } & (1 \leq k \leq n-2),\\
    \; & \; \\
    \displaystyle{\prod^{[\frac{n-1}{2}]}_{i=1}} 
     \dfrac{(q^{4i-2}pz_1/z_2;\xi^2)_{\infty} }{(q^{-4i}pz_1/z_2;\xi^2)_{\infty} } & (k=n-1), \\
    \; & \; \\
    \displaystyle{\prod^{[ \frac{n}{2} ]}_{i=1}}
     \dfrac{ (q^{4i-4}pz_1/z_2;\xi^2)_{\infty} }{ (q^{-4i+2}pz_1/z_2;\xi^2)_{\infty} } & (k=n).
  \end{cases}
\end{equation*}
Furthermore, $u_n^{(1)}$, $u_n^{(2)}$ and $ w^{(n,k)}(z)$ are given as follows:
\begin{eqnarray*}
u_n^{(1)} &=& (-q)^{ \frac{1}{2}n(n-1) } v_n,\\
u_n^{(2)} &=& (-q)^{ \frac{1}{2}(n-1)(n-2) } v_{++} \otimes v_{n-1},\\
\end{eqnarray*}
where $v_n$ is inductively defined by
\[
v_n = v_{+-}\otimes v_{n-1} + (-q)^{-n+1}v_{-+}\otimes {\sigma}_{n-1} v_{n-1}, \quad
v_1 = v_{+-}.
\]
\begin{equation*}
  w^{(n,k)}(z) = \sum^{ [ \frac{k}{2} ] }_{m=0} c_m^{(n,k)}(z) w_n^{(n-k-2m)},
\end{equation*}
\begin{eqnarray*}
& &c_m^{(n,k)}(z) 
  = \prod_{j=0}^{m-1}
      \frac{q^{-2(j+1)}}{(1+q^{2n-2k+4j})(1+q^{2n-2k+4j+2})(1-q^{2n-2k+2j})}\\
& & \qquad \qquad  \times  \sum_{i=0}^m \{ \prod_{j=0}^{m-1-i}(1+q^{2m-2j})
                                 \prod_{j=0}^{i-1}(1+q^{2n-2k+2m+2j})
                  {\fracwithdelims[][0pt]{m}{i}}_q((-1)^{n-k}q^{n+k-m}z)^i \},
\end{eqnarray*}
and $w_n^{(k)}$ is inductively defined by
\vspace{5mm}
\noindent for $n \geq 2$
\begin{eqnarray*}
  & w_n^{(k)}
    =\dfrac{-(-q)^{-n+k+2}}{(1+q^{2k})(1+q^{2k+2})}
          \{v_{+-+} \otimes w_{n-1}^{(k+1)}
          -(-q)^{-k-1}v_{-++}\otimes (\sigma_{n-1})_{12} \tilde{w}_{n-1}^{(k+1)} \}& \\
  & \hspace{-65mm} +v_{++-} \otimes w_{n-1}^{(k)}\quad & \hspace{-5mm}(1 \leq k \leq n),	\\
  & \quad \qquad +(-q)^{-n+k}\{v_{+-+} \otimes w_{n-1}^{(k-1)}
         +(-q)^{k-1}v_{-++}\otimes (\sigma_{n-1})_{12} \tilde{w}_{n-1}^{(k-1)} \}&\\
  & \; & \; \\
  & \tilde{w}_n^{(k)}
    =\dfrac{-(-q)^{-n+k+2}}{(1+q^{2k})(1+q^{2k+2})}
          \{v_{+-+} \otimes \tilde{w}_{n-1}^{(k+1)}
          -(-q)^{-k-1}v_{-++}\otimes (\sigma_{n-1})_{12} w_{n-1}^{(k+1)} \}& \\
  & \hspace{-65mm} +v_{++-} \otimes \tilde{w}_{n-1}^{(k)}	\quad &  \hspace{-5mm}(1 \leq k \leq n), \\
  & \quad \qquad +(-q)^{-n+k}\{v_{+-+} \otimes \tilde{w}_{n-1}^{(k-1)}
          +(-q)^{k-1}v_{-++}\otimes (\sigma_{n-1})_{12} w_{n-1}^{(k-1)}\}&\\
  & \; & \; \\
  & \hspace{-20mm} w_n^{(0)} 
    = \dfrac{-(-q)^{-n+2}}{(1+q^2)}
           \{v_{+-+} \otimes w_{n-1}^{(1)}
           -(-q)^{-1}v_{-++}\otimes (\sigma_{n-1})_{12} \tilde{w}_{n-1}^{(1)} \}& \\
  & \hspace{-65mm} +v_{++-} \otimes w_{n-1}^{(0)} \quad & (k = 0),	\\
  & \tilde{w}_n^{(0)}=0 \quad & (k = 0),\\
  & w_n^{(k)}=\tilde{w}_n^{(k)} = 0 & \hspace{-10mm} ( k > n \;\text{or}\; k < 0 ),
\end{eqnarray*}

\noindent for $n = 1$
\begin{eqnarray*}
 &w_1^{(k)} = \tilde{w}_1^{(k)} = v_{+-+}&( k = 1),\\
 &w_1^{(k)} = v_{++-},\;\;\tilde{w}_1^{(k)} = 0 &( k = 0), \\
 &w_1^{(k)} = \tilde{w}_1^{(k)} = 0  &(  k \not= 0,\;1).
\end{eqnarray*}
Here recursive formula of $w_n^{(k)}$ and $\tilde{w}_n^{(k)}$ are understood
as similar to Section 2.5 and 
$w_n^{(k)}$ denotes an equivalent class represented by itself in the space
\[
\{V^{(n)} \otimes V^{(n')} \otimes V^{(m)}\}/\{\operatorname{ker}T^{(k)} \otimes V^{(m)}\}
 = V^{(k)} \otimes V^{(m)},
\]
where $n' = \varphi^{(n)}(n-k)$ and $m = \varphi^{(n)}(k)$.

\section{Commutation relations}

We describe commutation relations of vertex operators 
and make a remark on difference equations for scalar functions which appeared in
the commutation relations. 

From the explicit forms of two point functions, 
we can obtain the following theorem. 
\begin{theorem}
  For any possible combination of weights $(\nu,\mu',\mu,\lambda)$ 
 (see Table 1 and 2 in Section 4.3),
    \begin{eqnarray}
      \Phi^{\nu \;V^{(n)}_2}_{\mu}(z_2) \Phi^{\mu \;V^{(n)}_1}_{\lambda}(z_1) 
        &=& r_n(z_1/z_2) \bar{R}^{(+,+)}_n(z_1/z_2)
          \Phi^{\nu \;V^{(n)}_1}_{\mu}(z_1) \Phi^{\mu \;V^{(n)}_2}_{\lambda}(z_2),
          \label{commutation relation 1} \\
      \Phi^{\nu \;V^{(n)}_2}_{\mu'}(z_2) \Phi^{V^{(k)}_1 \; \mu'}_{\lambda}(z_1)
        &=& \tau^{(k)}(z_1/z_2) 
          \Phi^{V^{(k)}_1 \; \nu}_{\mu}(z_1) \Phi^{\mu \;V^{(n)}_2}_{\lambda}(z_2),
          \label{commutation relation 2} \\
      \Phi^{\nu \;V^{(n)}_2}_{\mu'}(z_2) \Phi^{V^{(k)*a^{-1}}_1 \; \mu'}_{\lambda}(z_1)
        &=& \tau^{(k)}(z_1/z_2)^{-1}
          \Phi^{V^{(k)*a^{-1}}_1 \; \nu}_{\mu}(z_1) \Phi^{\mu \;V^{(n)}_2}_{\lambda}(z_2),
      \label{commutation relation 3}
    \end{eqnarray}
  where $r_n(z)$ and $\tau^{(k)}(z)$ are given by
    \begin{eqnarray*}
       && r_n(z) = z^{-\frac{n}{4}} \prod_{i=1}^{[ \frac{n}{2} ]}
         \dfrac{(q^{4i-2}z;\xi^2)_{\infty}(q^{4n-4i}z^{-1};\xi^2)_{\infty}}{
               (q^{4i-2}z^{-1};\xi^2)_{\infty}(q^{4n-4i}z;\xi^2)_{\infty}}, \\
       && \; \\ 
       && \tau^{(k)}(z) = 
        \begin{cases}
        {\displaystyle z^{-\frac{k}{2}} \prod_{j=1}^{k} }
         \dfrac{\Theta_{{\xi}^2}(-(-q)^{k+n-2j}z)}{
         \Theta_{{\xi}^2}(-(-q)^{k+n-2j}z^{-1})}
         \quad & (1 \leq k \leq n-2),\\ 
        \; & \; \\
        {\displaystyle z^{-\frac{n-2}{4}} \prod_{j=1}^{[ \frac{n-1}{2} ]} } 
         \dfrac{\Theta_{{\xi}^2}(-(-q)^{4j-1}z)}{\Theta_{{\xi}^2}(-(-q)^{4j-1}z^{-1})} 
         \quad & (k=n-1),\\ 
        \; & \; \\
        {\displaystyle  z^{-\frac{n}{4}} \prod_{j=1}^{[ \frac{n}{2} ]}  }
         \dfrac{\Theta_{{\xi}^2}(-(-q)^{4j-3}z)}{\Theta_{{\xi}^2}(-(-q)^{4j-3}z^{-1})}
         \quad & (k=n).
        \end{cases}
    \end{eqnarray*}
\end{theorem}
To prove this theorem, it is enough to show that  
the vacuum expectation values of the both sides in commutation relations are exactly same. 
Indeed, by irreducibility of $V(\lambda)$ and $V(\nu)$,
if the vacuum expectation values of both sides coincide, 
then the equality as operators on $V(\lambda)$ can be proved. 

This theorem gives us the explicit forms of the excitation spectra mentioned in Introduction.

We make a remark on a relation between 
difference equations for the scalar function $\tau^{(k)}(z)$
and fusion procedure of finite dimensional $U_q(D_n^{(1)})$-modules (cf. \cite{N}).

\begin{remark}
The scalar factor $\tau^{(k)}(z)$ satisfies the following 
difference equations:

\noindent for $n+k$:even
\begin{equation*}
  \begin{split}
   &\tau^{(k)}(z) 
     = (-1)^{-\frac{n}{2}}\tau^{(n)}(\zeta^{-1} z) \tau^{(n)}(\zeta z),  \\
   &\tau^{(k)}(z)  
     = (-1)^{-\frac{n}{2}+1}\tau^{(n-1)}(\zeta^{-1} z) \tau^{(n-1)}(\zeta z), 
  \end{split}
\end{equation*}

\noindent for $n+k$:odd
\begin{equation*}
  \begin{split}
    \tau^{(k)}(z)
      = (-1)^{-\frac{1}{2}(n-k-1)}\tau^{(n)}(\zeta^{-1} z)\tau^{(n-1)}(\zeta z),\\   
    \tau^{(k)}(z) 
      = (-1)^{-\frac{1}{2}(n-k-1)}\tau^{(n-1)}(\zeta^{-1} z) \tau^{(n)}(\zeta z),   
  \end{split}
\end{equation*}
where we put $\zeta = (-q)^{n-k-1}$.
We can see correspondence
between these equations and the following fusion procedure
\begin{eqnarray*}
    & V^{(k)}_z \hookrightarrow
    V^{(n)}_{\zeta^{-1} z} \otimes V^{(n)}_{\zeta z} \simeq
    V^{(n-1)}_{\zeta^{-1} z} \otimes V^{(n-1)}_{\zeta z}&  \quad ( n+k : \text{even} ),\\
    &  V^{(k)}_z \hookrightarrow
    V^{(n)}_{\zeta^{-1} z} \otimes V^{(n-1)}_{\zeta z} \simeq
    V^{(n-1)}_{\zeta^{-1} z} \otimes V^{(n)}_{\zeta z}&  \quad ( n+k : \text{odd} ).
\end{eqnarray*}
\end{remark}

\section{Formulation of the model}

We will construct the transfer matrix, the creation and annihilation operators, 
in the same way as in \cite{DFJMN}.

In quantum symmetry approach, we define the space of states by
\begin{eqnarray*}
{\cal F} = \underset{\lambda, \mu \in \{\Lambda_0, \Lambda_1, \Lambda_{n-1}, \Lambda_n \}}{
  \bigoplus} {\cal F}_{\lambda,\mu},
\end{eqnarray*}
where
\[
 {\cal F}_{\lambda,\mu} := V(\lambda) \otimes V(\mu)^{*a} 
    \equiv \operatorname{Hom}_{{\Bbb Q}(q)}(V(\mu), V(\lambda)).
\]
Following [1], we complete the space ${\cal F}_{\lambda,\mu}$ 
in the topology of formal power series in $q$. From now on we  
denote the completed space by the same simbol ${\cal F}_{\lambda,\mu}$.

The left $U_q(\frak g)$-module structure of ${\cal F}_{\lambda,\mu}$ can be written as
\[
x f = \sum x_{(1)} \circ f \circ a(x_{(2)}) \quad 
(f \in {\cal F}_{\lambda,\mu},\; x \in U_q(\frak g) ),
\]
where $\bigtriangleup(x) = \sum x_{(1)} \otimes x_{(2)}$.
At the same time, we can define right module structure on ${\cal F}_{\lambda,\mu}$ by
\[
f x = \sum a^{-1}(x_{(2)}) \circ f \circ x_{(1)},
\]
and then we denote this right module by ${\cal F}^{\; r}_{\lambda,\mu}$.
We have a natural pairing (cf. \cite{IIJMNT})
\begin{equation}
\langle f | g \rangle 
:= \frac{\operatorname{tr}_{V(\lambda)}(q^{-2\rho} f g )}{
         \operatorname{tr}_{V(\lambda)}(q^{-2\rho})}, \quad 
f \in {\cal F}_{\lambda,\mu}^{\; r}, \; g \in {\cal F}_{\mu,\lambda}.
\label{inner product}
\end{equation}
which satisfies
\[
\langle f x| g \rangle = \langle f | x g \rangle.
\]

In order to define ``local operators" on our space of states, 
we consider vertex operator
\begin{eqnarray*}
  \tilde{\Phi}_{\lambda \; V^{(n)}}^{\mu}(z) \; 
    : \; V(\lambda) \longrightarrow \hat{V}(\mu) \otimes V_z^{(k)}.
\end{eqnarray*}
We fix the normalization of $\tilde{\Phi}_{\lambda \; V^{(n)}}^{\mu}(z)$ by
\begin{eqnarray*}
  \tilde{\Phi}_{\Lambda_n \; V^{(n)}}^{\Lambda_0}(z) (|\Lambda_n \rangle \otimes v_{(+,+, \cdots,+)})
    = |\Lambda_0 \rangle + \cdots,
\end{eqnarray*}
and normalize the other vertex operators by using Dynkin diagram automorphisms.

Then we can show the following formula (cf. \cite{IIJMNT})
\begin{eqnarray*}
  \tilde{\Phi}_{\lambda \; V^{(n)}}^{\mu}(z) \tilde{\Phi}_{\lambda}^{\mu \; V^{(n)}}(z)
     &=& g \times \operatorname{id}_{V(\lambda)},\\
  \tilde{\Phi}_{\mu}^{\lambda \; V^{(n)}}(z) \tilde{\Phi}_{\mu \; V^{(n)}}^{\lambda}(z)
     &=& g \times \operatorname{id}_{V(\lambda) \otimes V^{(n)}}.
\end{eqnarray*}
where
\begin{eqnarray*}
  g =
  \begin{cases}
    {\displaystyle \prod_{i=1}^{\frac{n}{2}} }
    \dfrac{(q^{4i-2} \xi ; \xi^2)_{\infty}}{(q^{-4i+4} \xi ; \xi^2)_{\infty}}& \quad (n : \text{even}),\\
    \; & \; \\
    {\displaystyle \prod_{i=1}^{\frac{n-1}{2}} } 
    \dfrac{(q^{4i} \xi ; \xi^2)_{\infty}}{(q^{-4i+2} \xi ; \xi^2)_{\infty}}& \quad (n : \text{odd}).  
  \end{cases}                
\end{eqnarray*}
These formulas can be proved by the explicit forms of two point functions in Section 4.4 
in the same way as \cite{IIJMNT}.

By the formulas, we see that
\[
\tilde{\Phi}_{\lambda}^{\mu \; V^{(n)}}(z)\;:\;V(\lambda) \longrightarrow \hat{V}(\mu) \otimes V_z^{(n)}
\]
is invertible (cf. \cite{DFJMN}). 

We recall that, for any integral weight $\lambda$ of level one,
there is a unique level one weight $\mu$ such that non-trivial vertex operator
$\tilde{\Phi}^{\mu \; V^{(k)}}_{\lambda}(z)$ exists (see Section 4.1). We denote  such $\mu$
by $\lambda^{(k)}$. For vertex operator $\tilde{\Phi}^{\mu \; V^{(k)*a^{\pm 1}}}_{\lambda}(z)$,
we define $\lambda^{(k)*}$ similarly.

The row transfer matrix
\[
T_{\lambda,\mu}^{\lambda^{(n)},\mu^{(n)}} \; 
: \; {\cal F}_{\lambda,\mu} \longrightarrow {\cal F}_{\lambda^{(n)},\mu^{(n)}},
\]
is defined by the composition of the following operators
\begin{equation}
  V(\lambda) \otimes V(\mu)^{*a} 
  \longrightarrow V(\lambda^{(n)}) \otimes V^{(n)}_z \otimes V(\mu)^{*a} 
  \longrightarrow V(\lambda^{(n)}) \otimes V(\mu^{(n)})^{*a},
\end{equation}
where the first map is 
$\tilde{\Phi}^{\lambda^{(n)} \; V^{(n)}}_{\lambda}(z) \otimes \operatorname{id}$ and the second
one is $\operatorname{id} \otimes (\tilde{\Phi}^{\mu \; V^{(n)*a^{-1}}}_{\mu^{(n)}}(z))^t$.

Since $\operatorname{id}_{V(\lambda)} \in {\cal F}_{\lambda,\lambda}$, 
we define $| \text{vac} \rangle_{\lambda} := \operatorname{id}_{V(\lambda)}$.
In a similar way to \cite{IIJMNT}, we can show
\begin{eqnarray*}
T_{\lambda,\lambda}^{\lambda^{(n)},\lambda^{(n)}} | \text{vac} \rangle_{\lambda}
= g | \text{vac} \rangle_{\lambda^{(n)}},
\end{eqnarray*}  
from the next formula
\begin{eqnarray*}
 \bar{p}(\langle \tilde{\Phi}_{\mu}^{\lambda V^{(n)}_1}(z)
            \tilde{\Phi}_{\lambda}^{\mu V^{(n)*a^{-1}}_2}(z) \rangle ) &=& g,
\end{eqnarray*}  
where 
\begin{eqnarray*}
\bar{p}\;:\;V^{(n)} \otimes V^{(n)*a^{-1}} &\longrightarrow& {\Bbb Q}(q),\\
\quad v_1 \otimes v_2^* &\mapsto& \langle v_1 , v_2^* \rangle.
\end{eqnarray*}  

The creation and annihilation operators are constructed by using type II 
vertex operators. We express vertex operators as follows:
\begin{eqnarray*}
\Phi^{V^{(k)} \; \lambda^{(k)}}_{\lambda}(z)
&=& \sum_I 
v_{I} \otimes \Phi^{(k)}_{\lambda,I}(z), \\
\Phi^{V^{(k)*a^{-1}} \; \lambda^{(k)*}}_{\lambda}(z)
&=& \sum_I 
v_{I}^* \otimes \Phi^{(k)*}_{\lambda,I}(z), 
\end{eqnarray*}
where $\{v_I\}$ and $\{v_I^*\}$ are dual basis of 
$V^{(k)}$ and $V^{(k)*a^{-1}}$ and 
$\Phi^{(k)}_{\lambda,I}(z),\;
\Phi^{(k)*}_{\lambda,I}(z) \in {\cal F}_{\lambda,\lambda^{(k)*}}$.
The creation operator
\[
\phi^{(k)*}_{\lambda,I}(z) \; : \;
{\cal F}_{\lambda,\mu} \longrightarrow {\cal F}_{\lambda^{(k)*},\mu},
\]
is given by
\[
f \; \mapsto \Phi^{(k)*}_{\lambda,I}(z) \circ f.
\]
The annihilation operator 
\[
\phi^{(k)}_{\lambda,I}(z) \; : \;
{\cal F}_{\lambda,\mu} \longrightarrow {\cal F}_{\lambda^{(k)*},\mu},
\]
is defined by the adjoint of
\[
{\cal F}_{\mu, \lambda^{(k)}}^{\;r} \longrightarrow {\cal F}_{\mu, \lambda}^{\;r}, \quad
f \; \mapsto  f \circ \Phi^{(k)}_{\lambda,I}(z),
\]
with respect to the pairing \eqref{inner product}.

Commutation relations of vertex operators \eqref{commutation relation 2} and
\eqref{commutation relation 3} lead 
the following relations between the transfer matrix and 
the creation and annihilation operators:
\begin{eqnarray*}
\phi^{(k)*}_{\lambda^{(n)*},I}(z) T_{\lambda,\mu}^{\lambda^{(n)},\mu^{(n)}}
&=& \tau^{(k)}(z_1/z_2) 
T_{\lambda^{(k)*},\mu}^{(\lambda^{(k)*})^{(n)},\mu^{(n)}}\phi^{(k)*}_{\lambda,I}(z),\\
\phi^{(k)*}_{\lambda^{(n)},I}(z) T_{\lambda,\mu}^{\lambda^{(n)},\mu^{(n)}}
&=& \tau^{(k)}(z_1/z_2)^{-1} 
T_{\lambda^{(k)*},\mu}^{(\lambda^{(k)*})^{(n)},\mu^{(n)}}\phi^{(k)}_{\lambda,I}(z),
\end{eqnarray*}
where $\tau^{(k)}(z_1/z_2)$ is given in the previous section.

\section{Appendix A}

We calculate the scalar function
${\alpha}^{(\varepsilon_1,\varepsilon_2)}_n(z)$ in \eqref{second_inversion_spin}
and ${\beta}^{(\varepsilon_1,\varepsilon_2)}_n(z)$ in \eqref{image_MR_spin}.
We only consider ${\alpha}^{(+,+)}_n(z)$ and ${\beta}^{(+,+)}_n(z)$ for any even integer $n$.
(The other cases can be calculated similarly.) 

\subsection{Scalar function ${\alpha}^{(+,+)}_n(z)$}
The second inversion relation can be written as follows:
\begin{equation}
  {\alpha}^{(+,+)}_n(z)^{-1} (\bar{R}^{(+,+)}_n(z)^{-1})^{t_1}
   = (q^{2\rho}\otimes 1) (\bar{R}^{(+,+)}_n({\xi}^{-2}z)^{t_1})^{-1}
     (q^{-2\rho}\otimes 1).
  \label{second inversion 2}
\end{equation}
Let $v_n$ be a vector $v_{(+,\cdots,+)} \otimes v_{(-,\cdots,-)} \in V^{(+)} \otimes V^{(+)}$.
Applying both sides in \eqref{second inversion 2} to $v_n$, we can obtain
the explicit form of ${\alpha}^{(+,+)}_n(z)$.
We will show that
$v_n$ is an eigenvector for operators 
$(q^{2\rho}\otimes 1) (\bar{R}^{(+,+)}_n({\xi}^{-2}z)^{t_1})^{-1} (q^{-2\rho}\otimes 1)$
(resp. $(\bar{R}^{(+,+)}_n(z)^{-1})^{t_1}$)
and its eigenvalue is given by $\prod_{i=1}^{\frac{n}{2}} a(q^{4i-4}{\xi}^{-2}z)$
(resp. $\prod_{i=1}^{\frac{n}{2}} a(q^{4-4i}z^{-1})$).

By using the isomorphism \eqref{new_expression}, we see $v_n = v_{+-} \otimes v_{n-1}$ 
where $v_{n-1} \in V^{(+)}_{n-1} \otimes V^{(-)}_{n-1}$.
Taking transpose in the first component of the tensor product
in the recursive formulae \eqref{recursive form ++}, we can obtain recursive relations
for $\bar{R}^{(+,+)}_n(z)^{t_1}$. In particular we have
\begin{eqnarray}
\bar{R}^{(+,+)}_n(z)^{t_1} v_n
&=&a(z) v_{+-} \otimes \bar{R}^{(+,-)}_{n-1}(q^2z)^{t_1} v_{n-1},
\label{relation_App1}\\
\bar{R}^{(+,-)}_n(z)^{t_1} v_n
&=& v_{+-} \otimes \bar{R}^{(+,+)}_{n-1}(q^2z)^{t_1} v_{n-1},
\label{relation_App2}
\end{eqnarray}
where $a(z) = q(1-z)/(1-q^2z)$.
Combining \eqref{relation_App1} with \eqref{relation_App2}, we find
\[
\bar{R}^{(+,+)}_n(z)^{t_1} v_n 
=a(z) v_{+-} \otimes v_{+-} \otimes \bar{R}^{(+,+)}_{n-2}(q^4z)^{t_1} v_{n-2},
\]
then 
\begin{eqnarray*}
\bar{R}^{(+,+)}_n(z)^{t_1}(v_n) 
&=& \prod_{i=1}^{\frac{n}{2}} a(q^{4i-4}z) v_{+-} \otimes v_{+-} \otimes \cdots \otimes v_{+-}\\
&=& \prod_{i=1}^{\frac{n}{2}} a(q^{4i-4}z) v_n.
\end{eqnarray*}
Since $(q^{-2\bar{\rho}} \otimes 1)v_n = q^{-n(n-1)/2}v_n$,
we obtain 
\[
(q^{2\rho}\otimes 1) (\bar{R}^{(+,+)}_n({\xi}^{-2}z)^{t_1})^{-1}
  (q^{-2\rho}\otimes 1) v_n = \prod_{i=1}^{\frac{n}{2}} a(q^{4i-4}{\xi}^{-2}z)^{-1}v_n.
\]

We can similarly show
\[
(\bar{R}^{(+,+)}_n(z)^{-1})^{t_1} v_n = \prod_{i=1}^{\frac{n}{2}} a(q^{4i-4}z^{-1}) v_n,
\]
by using the first inversion relation 
\[
\bar{R}^{(+,+)}_n(z)^{-1} = P \bar{R}^{(+,+)}_n(z^{-1}) P \quad (P:\text{transposition}).
\]

Applying both sides in \eqref{second inversion 2} to $v_n$, we obtain
\begin{equation*}
  \begin{split}
    \alpha^{(+,+)}_n(z)
              &= \prod_{i=1}^{\frac{n}{2}} a(q^{4i-4}z^{-1}) a({\xi}^{-2}q^{4i-4}z)\\
              &={\displaystyle \prod_{i=1}^{\frac{n}{2}} }
                 \dfrac{(1-q^{-4i+4}z)(1-q^{-4n+4i}z)}{(1-q^{-4i+2}z)(1-q^{-4n+4i+2}z)}.
  \end{split}
\end{equation*}

\subsection{Scalar factor $\beta^{(+,+)}_n(z)$}
We can determine the explicit form of the scalar factor $\beta^{(+,+)}_n(z)$ in 
\eqref{image_MR_spin} by using $\alpha^{(+,+)}_n(z)$.

We recall that
\begin{equation}
  (\pi^{(+)}_{z_1} \otimes \pi^{(+)}_{z_2})({\cal R}'(z))
  = \beta^{(+,+)}_n(z) \bar{R}^{(+,+)}_n(z). 
 \label{R-mat(App)}
\end{equation}
and the factor $\beta^{(+,+)}_n(z)$ is a solution of 
the following difference equation (see Section 2.3):
\[
\alpha^{(+,+)}_n(z) \beta^{(+,+)}_n(z)^{-1} = \beta^{(+,+)}_n({\xi}^{-2}z)^{-1}.
\]
Then $\beta^{(+,+)}_n(z)$ can be written as 
\[
\beta^{(+,+)}_n(z) = c \prod_{i=1}^{\infty} \alpha^{(+,+)}_n(z \xi^{2i})^{-1} 
         = c {\displaystyle \prod_{i=1}^{ \frac{n}{2} } }
             \dfrac{(q^{4n-4i-2}z;{\xi}^2)_{\infty}(q^{4i-2}z;{\xi}^2)_{\infty}}{
             (q^{4n-4i}z;{\xi}^2)_{\infty}(q^{4i-4}z;{\xi}^2)_{\infty}}.
\]
for some constant $c$.
Applying the both sides in \eqref{R-mat(App)} to $v_{(+,\cdots,+)} \otimes v_{(+,\cdots,+)}$
we can determine the constant $c$.
Indeed, from the explicit form of the modified universal R-matrix in \eqref{modified R-matrix},
we find
\begin{equation*}
  \begin{split}
   &(\pi^{(+)}_{z_1} \otimes \pi^{(+)}_{z_2})({\cal R}'(z))
     (v_{(+,\cdots,+)} \otimes v_{(+,\cdots,+)})  \\
   &\quad = q^{-\sum_{i=1}^n \pi^{(+)} (h_i) \otimes \pi^{(+)} (\bar{\Lambda}_i) }
           (v_{(+,\cdots,+)} \otimes v_{(+,\cdots,+)}) + (\text{higher degree of $z$}),
  \end{split}
\end{equation*}
and by direct calculation we have
\begin{equation*}
  q^{-\sum_{i=1}^n \pi^{(+)}(h_i) \otimes \pi^{(+)}(\bar{\Lambda}_i)}
     (v_{(+,\cdots,+)} \otimes v_{(+,\cdots,+)}) 
    = q^{-\frac{n}{4}}(v_{(+,\cdots,+)} \otimes v_{(+,\cdots,+)}).
\end{equation*}
On the other hand we obtain
\begin{eqnarray*}
\beta^{(+,+)}_n(z) R^{(+,+)}_n(z) v_{(+,\cdots,+)} \otimes v_{(+,\cdots,+)}
&=& \beta^{(+,+)}_n(z) v_{(+,\cdots,+)} \otimes v_{(+,\cdots,+)}\\
&=& c v_{(+,\cdots,+)} \otimes v_{(+,\cdots,+)} + (\text{higher degree of $z$}),
\end{eqnarray*}
then we see
\[
c = q^{-\frac{n}{4}}.
\]

\section{Appendix B}

We describe calculation of two point function of the case 1 in Section 4.3.
(The other cases can be calculated similarly.)
In this case, $n$ is an even integer and two point function is given by
\begin{eqnarray*} 
  \langle \Phi^{ {\Lambda}_0 \;(V^{(n)})_2 }_{ {\Lambda}_n }(z_2)
   \Phi^{ {\Lambda}_n \; (V^{(n)})_1 }_{ {\Lambda}_0 }(z_1) \rangle.
\end{eqnarray*}
We denote this two point function by $\Psi_n(z_1,z_2)$. Here we remark that (cf. \cite{DFJMN})
\begin{eqnarray}
  \Psi_n(z_1,z_2)  \in V^{(+)} \otimes V^{(+)} \otimes (z_1/z_2)^{\Delta_{\Lambda_n}-\Delta_{\Lambda_0}}
  \otimes {\Bbb Q}(q)[[z_1/z_2]].
  \label{2-pt fun}
\end{eqnarray}

We can determine the explicit form of $\Psi_n(z_1,z_2)$ by using
Proposition 4.3 and the following lemma.
\begin{lemma}
If a vector $v_n$ in $V^{(+)} \otimes V^{(+)}$
satisfies 
\[
(\pi_{z_1}^{(+)} \otimes \pi_{z_2}^{(+)})
\bigtriangleup'(e_i)^{\langle h_i, \nu \rangle+1} v_n = 0,\;\; \text{and} \;\;
\operatorname{wt}v_n = 0,
\]
for all $i=0,1,\cdots,n$, then $v_n$ is uniquely determined up to multiple constant.
More exactly $v_n$ is given by 
\[
v_n = v_{+-}\otimes v_{n-1} + (-q)^{-n+1}v_{-+}\otimes {\sigma}_{n-1} v_{n-1}, \quad
v_1 = v_{+-}.
\]
\end{lemma}
\vspace{5mm}
Combining this lemma with Proposition 4.3, we can find 
\begin{equation}
\Psi_n(z_1,z_2) = \psi_n(z_1,z_2) v_n,
\label{solution}
\end{equation}
for some scalar function $\psi_n(z_1,z_2)$.
As being described in \cite{IIJMNT}, the two point function can be written as
\[
\Psi_n(z_1,z_2) = 
   z_1^{\Delta_{\Lambda_n}-\Delta_{\Lambda_0}}
   z_2^{\Delta_{\Lambda_0}-\Delta_{\Lambda_n}} \tilde{\Psi}_n(z_1/z_2).
\]
Then the scalar function $\psi(z_1,z_2)$ can be also written as
\begin{equation*}
  \psi_n(z_1,z_2) =
     z_1^{\Delta_{\Lambda_n}-\Delta_{\Lambda_0}}
     z_2^{\Delta_{\Lambda_0}-\Delta_{\Lambda_n}}
     \tilde{\psi}_n(z_1/z_2),
\end{equation*}
for some scalar function $\tilde{\psi}_n(z)$.
Substituting \eqref{solution} to the q-KZ equation we have
\begin{equation}
  (pz)^{\Delta_{\Lambda_n}} \tilde{\psi}_n(p z) v_n
    = z^{\Delta_{\Lambda_n}} \tilde{\psi}_n(z)  
      \beta^{(+,+)}_n(p z)\bar{R}^{(+,+)}_n(pz)(q^{-2\bar{\rho}} \otimes 1) v_n,
  \label{q-KZ 2}
\end{equation}
where $z = z_1/z_2$.

\begin{lemma} 
\begin{equation*}
   {\bar{R}^{(+,+)}_n(z)(q^{-2\bar{\rho}} \otimes 1) v_n = f_n(z) v_n }, 
\end{equation*}
\begin{eqnarray*}
f_n(z) =
 {\displaystyle q^{\frac{1}{2}n^2} \prod_{i=1}^{\frac{n}{2}} }
  \dfrac{(1-q^{-4i+2}z)}{(1-q^{4i-2}z)}.
\end{eqnarray*}
\end{lemma}
\noindent
Proof. \ 
Using the recursive relation 
\[
v_n = v_{+-}\otimes v_{n-1} + (-q)^{-n+1}v_{-+}\otimes {\sigma}_{n-1} v_{n-1}, 
\]
we can easily show the following formula:
\begin{equation}
X_{n-1}^{(+,-)}(q^{-2\bar{\rho'}} \otimes 1)v_{n-1} = [n-1]_qv_{n-1},
\label{App_formula}
\end{equation}
where $X_{n-1}^{(+,-)}$ is given by \eqref{def of X} and
$\bar{\rho'}=(2n-4)\omega_2+(2n-6)\omega_3+\cdots+2\omega_{n-1}$ (see Section 2.1).
By using the recursive relation of the R-matrices in \eqref{recursive form ++},
we find
\begin{equation*}
  \begin{split}
  \bar{R}^{(+,+)}_n(z)(q^{-2\bar{\rho}} \otimes 1) v_n
    =& q^{-n+1} v_{+-} \otimes a(z) \bar{R}^{(+,-)}_{n-1}(q^2 z) (q^{-2\bar{\rho'}} \otimes 1) v_{n-1} \\
     & + q^{-n+1} v_{-+} \otimes zb(z) \bar{R}^{(-,+)}_{n-1}(q^2 z) \sigma_{n-1}X^{(+,-)}_{n-1}
        (q^{-2\bar{\rho'}} \otimes 1) v_{n-1} \\
     & + (-1)^{-n+1} v_{+-} \otimes b(z) \bar{R}^{(+,-)}_{n-1}(q^2 z) \sigma_{n-1} X^{(-,+)}_{n-1}
        (q^{-2\bar{\rho'}} \otimes 1)  \sigma_{n-1} v_{n-1}  \\
     & + (-1)^{-n+1} v_{-+} \otimes a(z) \bar{R}^{(-,+)}_{n-1}(q^2 z) (q^{-2\bar{\rho'}} \otimes 1) 
         \sigma_{n-1} v_{n-1}.
  \end{split} 
\end{equation*}
Combining the above formula \eqref{App_formula} with 
\[ \sigma_{n-1} \bar{R}^{(+,-)}_{n-1}(z) = \bar{R}^{(-,+)}_{n-1}(z) \sigma_{n-1},\quad
   \sigma_{n-1} X^{(+,-)}_{n-1} = X^{(-,+)}_{n-1} \sigma_{n-1},
\]
and
\[
\sigma_{n-1} (q^{-2\bar{\rho'}} \otimes 1)  = (q^{-2\bar{\rho'}} \otimes 1) \sigma_{n-1},
\]
we have
\begin{equation*}
  \begin{split}
  \bar{R}^{(+,+)}_n(z)(q^{-2\bar{\rho}} \otimes 1) v_n
  = & q^n \frac{1-q^{-2n+2}z}{1-q^2z} 
         \{ v_{+-} \otimes \bar{R}^{(+,-)}_{n-1}(q^2 z) 
          (q^{-2\bar{\rho'}} \otimes 1) v_{n-1}\\  
    & \quad + (-q)^{-n+1}v_{-+} \otimes \sigma_{n-1} \bar{R}^{(+,-)}_{n-1}(q^2 z) 
        (q^{-2\bar{\rho'}} \otimes 1) v_{n-1} \}.
  \end{split} 
\end{equation*}
On the other hand, 
\begin{eqnarray*}
  \bar{R}^{(+,-)}_{n}(z)(q^{-2\bar{\rho}} \otimes 1) v_n
  &=& q^{n-1} \frac{1-q^{-2n+2}z}{1-z} 
         \{ v_{+-} \otimes \bar{R}^{(+,+)}_{n-1}(q^2 z) 
          (q^{-2\bar{\rho'}} \otimes 1) v_{n-1}\\  
  & &\quad + (-q)^{-n+1}v_{-+} \otimes \sigma_{n-1} \bar{R}^{(+,+)}_{n-1}(q^2 z) 
        (q^{-2\bar{\rho'}} \otimes 1) v_{n-1} \},
\end{eqnarray*}
for any odd integer $n$.
We obtain the eigenvalue $f_n(z)$ by using these recursive relation inductively.
\hfill Q.E.D
\vspace{5mm}

By using this lemma, the q-KZ equation \eqref{q-KZ 2} can be reduced to the following 
difference equation:
\begin{equation*}
\tilde{\psi}_n(pz)
 = p^{-\Delta_{\Lambda_n}} \beta^{(+,+)}_n(p z) f_n(p z)\tilde{\psi}_n(z).
\end{equation*}
Then we have
\begin{equation*}
 \begin{split}
   \tilde{\psi}_n(z) 
    & = c \prod_{j=1}^{\infty}\{ p^{\Delta_{\Lambda_n}} 
               f_n(p^j z)^{-1}\beta_n(p^j z)^{-1} \} \\
    & = c \prod_{j=1}^{\infty} \{ 
               {\displaystyle \prod_{i=1}^{\frac{n}{2}} }
                  \dfrac{(1-q^{4i-2}p^jz)}{(1-q^{-4i+2}p^jz)}
               {\displaystyle \prod_{i=1}^{ \frac{n}{2} } }
                  \dfrac{(q^{4n-4i}p^jz;{\xi}^2)_{\infty}(q^{4i-4}p^jz;{\xi}^2)_{\infty}}{
                         (q^{4n-4i-2}p^jz;{\xi}^2)_{\infty}(q^{4i-2}p^jz;{\xi}^2)_{\infty}} 
                            \} \\
    & = c \prod_{j=1}^{\infty} \{ 
               {\displaystyle \prod_{i=1}^{ \frac{n}{2} } }
                  \dfrac{(q^{-4i+2}p^{j+1}z;{\xi}^2)_{\infty}(q^{4i-4}p^jz;{\xi}^2)_{\infty}}{
                         (q^{-4i+2}p^jz;{\xi}^2)_{\infty}(q^{4i-4}p^{j+1}z;{\xi}^2)_{\infty}} 
                            \} \\
    & = c \prod^{ \frac{n}{2} }_{i=1}
             \frac{ (q^{4i-4}pz;\xi^2)_{\infty} }{ (q^{-4i+2}pz;\xi^2)_{\infty} }.
  \end{split}
\end{equation*}
By the property in \eqref{2-pt fun}, we see that
\[
  \tilde{\psi}_n(z) \in {\Bbb Q}(q)[[z]].
\]
Then $c$ is a constant. This constant $c$ can be determined 
by the constant term of the two point function $\tilde{\Psi}_n(z_1/z_2) = \tilde{\psi}_n(z_1/z_2) v_n$.
It is easy to know this constant term 
from the normalization of vertex operators (see Section 4.1) and we find
\[
c = (-q)^{\frac{1}{2}n(n-1)}.
\]
Then we obtain the explicit form of the two point function.

\section{Appendix C}

We give explicit forms of isomorphisms of $U_q(\frak g)$-module in \eqref{dual for fusion}
\begin{equation*}
C_{\pm}^{(k)} \;
: \; V_{z\xi^{\mp}}^{(k)} \longrightarrow V_z^{(k)*a^{\pm 1}} \quad (1\leq k \leq n-2),
\end{equation*}
though the existence of these isomorphisms is proved in \cite{DaiO}.

Let us recall that there exist the following isomorphisms of $U_q(\overset{\circ}{\frak g})$-modules
\begin{equation}
V_z^{(k)} \simeq 
V_{\bar{\Lambda}_{k}} \oplus V_{\bar{\Lambda}_{k-2}} \oplus
\cdots \oplus (V_{\bar{\Lambda}_{1}} \; \text{or} \; V_{0})\;\;(k = 1,2,\ldots n-2),
\label{dual_iso}
\end{equation}
where $V_{\bar{\Lambda}_{i}}$ is denoted the irreducible $U_q(\overset{\circ}{\frak g})$-module
with highest weight $\bar{\Lambda}_{i}$ (cf. \cite{O}).
Here we know that 
if two irreducible finite-dimensional representations of $U_q(\overset{\circ}{\frak g})$ 
have the same highest weight
then they are isomorphic. 
For $i = 1,2,\ldots n-2$, $V_{\bar{\Lambda}_{i}}$ and 
$(V_{\bar{\Lambda}_{i}})^{*a^{\pm1}}$ have the same highest weight, 
then there exists an isomorphism
\[
\bar{C}_{i,\pm}:V_{\bar{\Lambda}_{i}} \longrightarrow (V_{\bar{\Lambda}_{i}})^{*a^{\pm1}},
\]
and this isomorphism is unique up to multiple constant.
Thanks to the isomorphism \eqref{dual_iso},
we can write $C_{\pm}^{(k)}$ as a linear combination of $\bar{C}_{i,\pm}$.
Then we normalize $\bar{C}_{i,\pm}$ and express the isomorphism $C_{\pm}^{(k)}$
by using them. 
Let us recursively define a vector $x_m$ by
\begin{equation*}
x_m = v_{+-} \otimes x_{m-1} + (-q)^{m-1}v_{-+} \otimes \sigma_{m-1} x_{m-1},
\quad x_1 = v_{+-},
\end{equation*}
where $x_m \in V^{(+)}_m \otimes V^{(\varepsilon)}_m$ and if $m$ is even (resp. odd) then 
$\varepsilon = +,\;(\text{resp.}-)$. 
By using $x_m$, we recursively define $y_{n,\pm}^{(m)} \;(m \leq n)$ by
\begin{eqnarray*}
    &\hspace{-8mm} y_{n,+}^{(m)} = v_{++} \otimes y_{n-1,+}^{(m)} \quad & (n>m),\\
    &y_{n,-}^{(m)} = v_{--} \otimes \sigma_{n-1} y_{n-1,-}^{(m)} \quad & (n>m),\\
    &y_{m,+}^{(m)} = y_{m,-}^{(m)} = x_{m} \quad &(m=n).
\end{eqnarray*}
We remark that  
\begin{equation*}
V^{(k)} = V^{(n)} \otimes V^{(n')}/\operatorname{ker}T^{(k)} \quad (n' = \varphi(n-k) ),
\end{equation*}
and the equivalent classes represented by $y_{n,+}^{(m)}$ (resp. $y_{n,-}^{(m)}$) 
are the highest weight vector
(resp. the lowest weight vector) in $U_q(\overset{\circ}{\frak g})$-module
$V_{\bar{\Lambda}_{n-m}}$.
Let us normalize $\bar{C}_{i,\pm} \;(i = 1,2,\ldots,n-2)$ by
\begin{equation*}
\bar{C}_{n-m,\pm}(y_{m,+}) = y_{m,-}^* \;(m = 2,3,\ldots,n).
\end{equation*}
By using these isomorphisms $\bar{C}_{i,\pm} \;(i = 1,2,\ldots,n-2)$ we have
\begin{equation*}
C_{\pm}^{(k)} = \sum_{m=0}^{[ \frac{k}{2} ]}
               k_{m,\pm} \bar{C}_{k-2m,\pm}.
\end{equation*}
By direct calculation we can show 
\begin{eqnarray*}
  k_{m,\pm} &=& (-q)^{\pm m(2n-2k-2m+1)} \prod_{i=0}^{m-1} 
      \frac{[2n-2k+2i]_q}{[2i+2]_q} (1+q^{2n-2k+4i})(1+q^{2n-2k+4i+2}).
\end{eqnarray*}


\begin{thebibliography}{99}
\bibitem{DFJMN} B. Davies, O. Foda, M. Jimbo, T. Miwa, and A. Nakayashiki: 
Diagonalization of the XXZ Hamiltonian by vertex operators, \;
Comm. Math. Phys., {\bf 151}, (1993), 86--153.
\bibitem{IIJMNT} M. Idzumi, K. Iohara, M. Jimbo, T. Miwa, T. Nakashima and T. Tokihiro: 
Quantum affine symmetry in vertex models, \;
Int. J. Mod. Phys., A, {\bf 8}, (1993), 1479--1511.
\bibitem{DatO} E. Date and M. Okado: 
Calculation of excitation spectra of the spin model related with 
the vector representation of the quantized affine algebra of type $A_n^{(1)}$, \;
Int. J. Mod. Phys., A, {\bf 9}, (1994), 399--417.
\bibitem{DJO} E. Date M. Jimbo and M. Okado: 
Crystal base and q-vertex operators, \;
Comm. Math. Phys., {\bf 155}, (1993), 47--69.
\bibitem{DaiO} B. Davies and M. Okado: 
Excitation spectra of the spin models constructed from quantized affine algebras 
of type $B_n^{(1)},D_n^{(1)}$, \;
Int. J. Mod. Phys., A, {\bf 11}, (1996), 1975--2017.
\bibitem{O} M. Okado: 
Quantum $R$ matrix related to the spin representation of $B_n^{(1)},D_n^{(1)}$, \;
Comm. Math. Phys., {\bf 134}, (1990), 467--486.
\bibitem{KKMMNN} S-J. Kang, M. Kashiwara, K. C. Misra, T. Miwa, T. Nakashima and A. Nakayashiki:
Affine crystals and vertex models, Int. J. Mod. Phys. A, {\bf 7}, Suppl. 1A (1992)
449--484.
\bibitem{BVV} O. Babelon, H. J. de Vega and C. M. Viallet: 
Exact excitation spectra of the $Z_{n+1} \times Z_{n+1}$ generalized Heisenberg model,\;
Nucl. Phys. B, {\bf 220}, [FS8] (1983), 283--301.
\bibitem{N} T. Nakanishi: 
Fusion, mass and representation theory of the Yangian algebra,\;
Nucl. Phys. B, {\bf 439}, (1995), 441--459.
\bibitem{ORW} E. Ogievesky, N. Yu. Reshetikhin and P. Wiegmann:
The principal chiral field in two dimensions on classical Lie algebras,
Nucl. Phys. B, {\bf 15}, (1987), 45--98.
\bibitem{K} V. G. Kac:
Infinite dimensional Lie algebras,\;
3rd ed., Cambridge University Press, Cambridge 1990.
\end{thebibliography}
\end{document}